\documentclass[12pt,psamsfonts]{article}
\usepackage[dvips]{epsfig}
\usepackage{graphics}
\usepackage{amsmath}
\usepackage{amsfonts}
\usepackage{amssymb}
\usepackage{amsthm}
\usepackage{bm}
\usepackage{enumerate}
\usepackage[mathscr]{eucal}
\usepackage{inputenc}
\usepackage[T2A]{fontenc}
\theoremstyle{plain}
\newtheorem{theorem}{Theorem}[section]
\newtheorem{lemma}{Lemma}[section]
\newtheorem{definition}{Definition}[section]
\newtheorem{remark}{Remark}[section]

\newtheorem*{nonumtheorem}{Theorem}

 \numberwithin{equation}{section}

\begin{document}
\title{\textbf{The truncated Fourier operator. IV.\\
}} \footnotetext{\hspace*{-4.0ex}\textbf{Mathematics Subject
Classification: (2000).} Primary 47E05, 34E05; Secondary
33E10.\endgraf \hspace*{-2.0ex}\textbf{Keywords:} Truncated
Fourier operator, prolate spheroid differential operator,
selfadjoint extensions of  singulary differential operators,
abstract boundary conditions, commuting operators.}
\author{Victor Katsnelson \and Ronny Machluf}
\date{\ }%
 \maketitle
 \abstract{We consider the formal prolate spheroid differential
operator on a finite symmetric interval and describe all its
self-adjoint boundary conditions. Only one of these boundary
conditions corresponds to a self-operator differential operator
which commutes with the Fourier operator truncated on the
considered finite symmetric interval.}

\setcounter{section}{3}
\section{Self-adjoint boundary conditions~for \\the prolate spheroid differential
\mbox{operator.}}

The study of the spectral theory of the Fourier operator
restricted on a finite symmetric interval \([-a,a]\):
\begin{multline}
\label{TrFSI}%
 \hspace{4.0ex}
(\mathscr{F}_{E}x)(t)=\frac{1}{\sqrt{2\pi}}\int\limits_{-a}^{a}e^{it\xi}x(\xi)\,d\xi,\quad
\ \ t\in{}E,\ \ E=[-a,a],\\
\mathscr{F}_{E}:\,\,L^2(E)\to{}L^2(E)\,,\hspace{4.0ex}
\end{multline}
is closely related to study of the differential operator generated
by the differential expression (or \emph{formal differential operator}) \(L\):
\begin{equation}
\label{FDOFI}%
(Lx)(t)=
-\frac{d\,\,}{dt}\bigg(1-\frac{t^2}{a^2}\bigg)\frac{dx(t)}{dt}+t^2x(t)\,.
\end{equation}
The operator \(L\) is said to be the \textsf{prolate spheroid
differential operator.}

The relationship between the spectral theory theories of the
integral operator
\(\mathscr{F}_{E}^{\ast}\mathscr{F}_{E},\,\,E=[-a,a],\) and the
prolate spheroid differential operator was discovered in the
series of remarkable papers \cite{SlPo}, \cite{LaP1}, \cite{LaP2},
where this relationship has been ingeniously used
 for developing the spectral theory of the operator
\(\mathscr{F}_{E}^{\ast}\mathscr{F}_{E}\). (See also \cite{Sl2},
\cite{Sl3}.) Actually the reasoning of  \cite{SlPo}, \cite{LaP1},
\cite{LaP2} can be easily applied to the spectral theory of the
operator \(\mathscr{F}_{E}\) itself rather the operator
\(\mathscr{F}_{E}^{\ast}\mathscr{F}_{E}\).

It should be emphasized that what was used in \cite{SlPo},
\cite{LaP1}, \cite{LaP2} this is a certain system of
eigenfunctions related to the differential expression \(L\),
\eqref{FDOFI}. These eigenfunctions are known as \emph{prolate
wave functions}. The prolate wave functions themselves were used
much before the series of the papers \cite{SlPo}, \cite{LaP1},
\cite{LaP2} was published. These functions naturally appears by
separation of variables in the Laplace equation in spheroidal
coordinates. However this was the work \cite{SlPo}, \cite{LaP1},
\cite{LaP2} where the prolate functions were first used for
solving  the spectral problem related to the Fourier analysis on a
finite symmetric interval. Until now, there is no clear
understanding why the approach used in \cite{SlPo}, \cite{LaP1},
\cite{LaP2} works. This is a a lucky accident which still waits
for its explanation. (See \cite{Sl3}.)

\emph{Actually eigenfunctions are related not to the the
differential expression itself but to a certain differential
operator generated by the differential expression. This
differential operator is generated not only by the differential
expression but also by certain boundary condition.} In the case
\(E=(-\infty,\infty)\), the differential operator generated by the
differential expression \(-\frac{d^2\,}{dt^2}+t^2\) on the class
smooth finite functions (or the class of smooth fast decaying
functions) is essentially selfadjoint: the closure of this
operator is a selfadjoint operator. Thus in the case
\(E=(-\infty,\infty)\) there is no need to discuss the boundary
condition because there is no such boundary conditions.

In contrast to the case \(E=(-\infty,\infty)\), in the case
\(E=[-a,a],\,\,0<a<\infty\) the minimal differential operator
\(\displaystyle-\frac{d\,\,}{dt}\bigg(1-\frac{t^2}{a^2}\bigg)\frac{dx\,\,}{dt}+t^2\)
is symmetric but is \emph{not} self-adjoint. This minimal operator
admits the family of self-adjoint extensions. Each of this
selfadjoint extensions is described by a certain boundary
conditions at the end points of the interval \([-a,a]\). The set
of all such extensions can be parameterized by the set of all
\(2\times2\) unitary matrices.

 It turns out that only one of these extensions
commutes with the truncated Fourier operator
\(\mathscr{F}_{E},\,\,E=[-a,a]\). To our best knowledge, until now
no attention was paid to this aspect. In the present paper, we in
particular investigate the question which extensions of the
minimal differential operator generated by
\(L\), \eqref{FDOFI}, commute with \(L\).\\[2.0ex]
\hspace*{2.0ex} \textsf{Analysis of solutions of the equation
\(Lx=\lambda{}x\) near singular points.}

 For the differential equation
\begin{equation}
\label{DEFEWP}
-\frac{d\,\,}{dt}\bigg(1-\frac{t^2}{a^2}\bigg)\frac{dx(t)}{dt}+t^2x(t)=\lambda{}x(t),\quad
 t\in\mathbb{C},
\end{equation}
considered in complex plane, the points \(-a\) and \(a\) are the
regular singular point. Let us investigate the asymptotic behavior
of solutions of this equation near these points. (Actually we need
to know this behavior only for real \(t\in(-a,a)\), but it is much
easier to investigate this question using some knowledge from the
analytic theory of differential equation.) Concerning the analytic
theory of differential equation see \cite[Chapter 5]{Sm}.

Let us outline an analysis of solution of the equation near the
point \(t=-a\). Change of variable
\[t=-a+s,\ \ \  x(-a+s)=y(s)\]
reduces the equation \eqref{DEFEWP} to the form
\begin{equation}
\label{ChVar}%
 s\frac{d^2y(s)}{ds^2}+p(s)\frac{dy(s)}{ds}+q(s,\lambda)y(s)=0\,,
\end{equation}
where \(p(s)\) and \(q(s)\) are functions holomorphic within the
disc \(|s|<2a\), moreover \(p(0)=1\):
\begin{equation}
\label{psex}%
p(s)=1+\sum\limits_{k=1}^{\infty}p_ks^k,\quad
q(s)=\sum\limits_{k=0}^{\infty}q_k(\lambda)s^k\,.
\end{equation}
An explicit calculation with power series give:
\begin{equation}
\label{ExEx1} p_1=-\frac{1}{2a};\ \ \
q_0=\frac{\lambda{}a}{2}-\frac{a^3}{2},\ \ q_1=
\frac{\lambda}{4}+\frac{3}{4}a^2\,.
\end{equation}
 Now we turn to the analytic theory of differential equations.
 The results of this theory which we need are
 presented for example in \cite[Chapter 5]{Sm}, see especially
section \textbf{98} there. We seek the solution of the equation
\eqref{ChVar}-\eqref{psex} in the form
\[y(s)=s^{\rho}\sum\limits_{k=0}^{\infty}c_ks^k\,.\]
Substituting this to the left-hand side of the equation
\eqref{ChVar}-\eqref{psex} and equating the coefficients of like
powers of \(s\) to zero we obtain the equations for the
determination of \(\rho\) and \(c_k\). In particular, the equation
corresponding to the power \(s^0\) is of the form:
\[c_0\,\rho^2=0\,.\]
The coefficient \(c_0\) plays the role of a normalizing  constant,
and we may take
\begin{equation}
\label{NoCo}%
 c_0=1\,.
\end{equation}
Equation for \(\rho\), the so called \emph{characteristic
equation}, is of the form
\begin{equation}
\label{ChaEq}%
 \rho^2=0.
\end{equation}
This equation has the root \(\rho=0\), and this root is multiple.
According to general theory, the equation
\eqref{ChVar}-\eqref{psex} has two solutions \(y_1(s)\) and
\(y_2(s)\) possessing the properties:

The solution \(y_1(s)\) is a function holomorphic is the disc
\(|s|<2a\) satisfying the normalizing condition \(y_1(0)=1\). The
solution \(y_2(s)\) is of the form \(y_2(s)=y_1(s)\,\ln{}s+z(s)\),
where \(z(s)\) is a function holomorphic  in the disc \(|s|<2a\)
and satisfying the condition \(z(0)=0\). WE may calculate
explicitly several first coefficients of power expansions
\[y_1(s)=1+\sum\limits_{k=1}^{\infty}c_ks^k,\ \ \ z(s)=\sum\limits_{k=1}^{\infty}d_ks^k\,:\]
\[c_1=\frac{a^3}{2}-\frac{\lambda{}a}{2},\ \ \
d_1=\lambda{}a-a^3+\frac{1}{2a}\,.\] Returning to the variable
\(t=-a+s\), we get the following result:
\begin{lemma}%
\label{ABSNS}%
Let \(L\) be the differential expression defined by \eqref{FDOFI},
and \(\lambda\in\mathbb{C}\) be arbitrary fixed.\\[1.0ex]
\hspace*{2.0ex}\textup{1.}%
\begin{minipage}[t]{0.93\linewidth}
There exist two solutions \(x_{1}^{-}(t,\lambda)\) and
\(x_{2}^{-}(t,\lambda)\) of
the equation %
\mbox{ \begin{math}%
Lx(t)=\lambda{}x(t)
\end{math}} %
possessing the properties:
\hspace*{1.0ex}\textup{a.}%
\begin{minipage}[t]{0.96\linewidth} The
function \(x_{1}^{-}(t,\lambda)\) is holomorphic in the disc
\(|t+a|<2a\), and satisfy the normalizing condition
\(x_{1}^{-}(-a,\lambda)=1\)\,;
\end{minipage}\\
\hspace*{1.0ex}\textup{b.}%
\begin{minipage}[t]{0.98\linewidth}
 The function
\(x_{2}^{-}(t,\lambda)\) is of the form %
\[x_{2}^{-}(t,\lambda)=x_{1}^{-}(t,\lambda)\,\ln{}(t+a)+w^{-}(t,\lambda),\]
where the function \(w^{-}(t,\lambda)\) is holomorphic in the disc
\(|t+a|<2a\) and satisfy the condition \(w^{-}(-a,\lambda)=0\)\,.
\end{minipage}
\end{minipage}\\
\hspace*{2.0ex}\textup{2.}
\begin{minipage}[t]{0.96\linewidth}
There exist two solutions \(x_{1}^{+}(t,\lambda)\) and
\(x_{2}^{+}(t,\lambda)\) of
the equation %
 \mbox{\begin{math}%
Lx(t)=\lambda{}x(t)
\end{math}} %
possessing the properties:
\hspace*{1.0ex}\textup{a.}%
\begin{minipage}[t]{0.96\linewidth} The
function \(x_{1}^{+}(t,\lambda)\) is holomorphic in the disc
\(|t-a|<2a\), and satisfy the normalizing condition
\(x_{1}^{+}(a,\lambda)=1\)\,;
\end{minipage}\\
\hspace*{1.0ex}\textup{b.}
\begin{minipage}[t]{0.96\linewidth}
 The function
\(x_{2}^{+}(t,\lambda)\) is of the form %
\[x_{2}^{+}(t,\lambda)=x_{1}^{+}(t,\lambda)\,\ln{}(a-t)+w^{+}(t,\lambda),\]
where the function \(w^{+}(t,\lambda)\) is holomorphic in the disc
\(|t+a|<2a\) and satisfy the condition \(w^{+}(a,\lambda)=0\)\,.
\end{minipage}
\end{minipage}%
\end{lemma}%

Given a fixed \(\lambda\), the solutions \(x_{1}^{-}(t,\lambda)\),
\(x_{2}^{-}(t,\lambda)\) are linearly independent, therefore
arbitrary solution \(x(t,\lambda)\) of the equation \eqref{DEFEWP}
can be expanded into a linear combination
\begin{subequations}
\label{LCom}
\begin{equation}
\label{LComM}%
x(t,\lambda)=c_1^{-}x_{1}^{-}(t,\lambda)+c_2^{-}x_{2}^{-}(t,\lambda).
\end{equation}
 The solutions \(x_{1}^{+}(t,\lambda)\),
\(x_{2}^{+}(t,\lambda)\) also are linearly independent, and the
solution \(x(t,\lambda)\) can be also expanded into the other
linear combination
\begin{equation}
\label{LComP}%
x(t,\lambda)=c_1^{+}x_{1}^{+}(t,\lambda)+c_2^{+}x_{2}^{+}(t,\lambda).
\end{equation}
\end{subequations}
Here \(c_1^{\pm}\), \(c_2^{\pm}\) are constants (with respect to
\(t\)). The solution \(x_{1}^{-}(t,\lambda)\) is bounded  and the
solution \(x_{2}^{-}(t,\lambda)\) grows logarithmically as
\(t\to{}-a\). Therefore the solution \(x(t,\lambda)\) is square
integrable near the point \(t=-a\). For the same reason, the the
solution \(x(t,\lambda)\) is square integrable near the point
\(t=a\). Thus we prove the following result.
\begin{lemma}\label{InArS}
Given \(\lambda\in\mathbb{C}\), then every solution
\(x(t,\lambda)\) of the equation \eqref{DEFEWP} satisfy the
condition
\begin{equation}
\label{SqInt}%
 \int\limits_{-a}^{a}\big|x(t,\lambda)\big|^2\,dt<\infty\,.
\end{equation}
\end{lemma}

\begin{center}
 \textsf{Differential operators related to the differential
expression \(L\), \eqref{DiffExp}.}
\end{center}

With the differential expression  (or, in other words, the \emph{formal  differential operator}) \(L\),
\begin{equation}
\label{DiffExp}%
L=-\frac{d\,\,}{dt}\bigg(1-\frac{t^2}{a^2}\bigg)\frac{d\,\,}{dt}+t^2\,,
\end{equation}
various differential operators may be related according to whether
boundary conditions are posed on functions from their domains of
definition.

\begin{definition}
\label{DoDeMaDO}%
 The set \(\mathcal{A}\) is the set of complex-valued
functions \(x(t)\) defined on the open interval \((-a,a)\) and
satisfied the
following conditions:\\
\hspace*{1.5ex}\textup{1.}\,
\begin{minipage}[t]{0.95\linewidth} The derivative
\(\dfrac{dx(t)}{dt}\) of the function \(x(t)\) exists at every
point \(t\) of the interval \((-a,a)\);
\end{minipage}\\
\hspace*{1.5ex}\textup{2.}\,
\begin{minipage}[t]{0.95\linewidth}
The function \(\dfrac{dx(t)}{dt}\) is absolutely continuous on
every compact subinterval of the interval \((-a,a)\);
\end{minipage}\\
\end{definition}
\begin{definition}
\label{DoDeMiDO}%
The set \(\mathring{\mathcal{A}}\)
is the set of complex-valued functions \(x(t)\) defined on the open interval
\((-a,a)\) and satisfied the
following conditions:\\
\hspace*{1.5ex}\textup{1.}\,
\begin{minipage}[t]{0.95\linewidth}
The function \(x(t)\) belongs to the set \(\mathcal{A}\) defined
above;
\end{minipage}\\[1.0ex]
\hspace*{1.5ex}\textup{2.}\,
\begin{minipage}[t]{0.95\linewidth}
The support \(\textup{supp}\,x\) of the function \(x(t)\) is a
compact subset of the open interval \((-a,a)\):
\((\textup{supp}\,x)\!\Subset{}(-a,a)\).
\end{minipage}
\end{definition}

\begin{definition}
\label{MaxDO}
The differential operator \(\mathcal{L}_{\textup{max}}\) is defined as follows:\\[1.0ex]
\hspace*{1.5ex}\textup{1.}\,
\begin{subequations}
\label{maxdo}
\begin{minipage}[t]{0.95\linewidth}
The domain of definition \(\mathcal{D}_{\mathcal{L}_{\textup{max}}}\)%
of the operator \(\mathcal{L}_{\textup{max}}\)
is:
\begin{equation}%
\label{maxdo1}
\mathcal{D}_{\mathcal{L}_{\textup{max}}}=\lbrace{}x:\,x(t)\in%
L^2((-a,a))\cap{\mathcal{A}}\ \  \textup{and} \ \
 (Lx)(t)\in{}L^2((-a,a))\rbrace,
 \end{equation}%
 where \((Lx)(t)\) is defined\,\footnotemark%
by \eqref{FDOFI}.
\end{minipage}\\[1.0ex]
\hspace*{1.5ex}\textup{2.}\,
\begin{minipage}[t]{0.95\linewidth}
The action of the operator  \(\mathcal{L}_{\textup{max}}\) is:
\begin{equation}
\label{maxdo2}
\textup{For}\ x\in\mathcal{D}_{\mathcal{L}_{\textup{max}}}\,,\,\ %
\mathcal{L}_{{}_\textup{max}}x=Lx\,.
\end{equation}
\end{minipage}{\ }\\[2.0ex]
\end{subequations}
\footnotetext{Since \(x\in\mathcal{A}\), the expression \((Lx)(t)\) is well defined.}%
The operator \(\mathcal{L}_{\textup{max}}\) is said to be the \emph{maximal differential
operator generated by the differential expression} \(L\), \eqref{DiffExp}.
\end{definition}
The minimal differential operator
\(\mathcal{L}_{{}_\textup{min}}\) is the restriction of the
maximal differential operator \(\mathcal{L}_{{}_\textup{max}}\) on
the set of functions which is some sense vanish at the endpoint of
the interval \((-a,a)\). The precise definition is presented
below.
\begin{definition}
\label{MinDO}
\begin{subequations}
\label{mindo}
The operator \(\mathcal{L}_{{}_\textup{min}}\)
is the closure\,\footnote{%
Since the operator \(\mathring{\mathcal{L}}\) is symmetric, it  is closable.}%
 of the operator \(\mathring{\mathcal{L}}\):
\begin{equation}
\label{mindo1}
\mathcal{L}_{{}_\textup{min}}=
\textup{clos}\big(\mathring{\mathcal{L}}\,\big)\,,
\end{equation}
where the operator \(\mathring{\mathcal{L}}\) is the restriction of the
operator \(\mathcal{L}_{{}_\textup{max}}\):
\begin{equation}%
\label{mindo2}
\mathring{\mathcal{L}}\subset{}\mathcal{L}_{\textup{max}},\quad %
\mathring{\mathcal{L}}=%
{\mathcal{L}_{\textup{max}}}_{|_{\scriptstyle\mathcal{D}_{\mathring{\mathcal{L}}}}},
\quad
\mathcal{D}_{\mathring{\mathcal{L}}}=%
\mathcal{D}_{\!{}_{\scriptstyle\mathcal{L}_{\textup{max}}}}\!\cap{}\mathring{\mathcal{A}}\,.
\end{equation}%
\end{subequations}
\end{definition}

By \(\langle\,,\,\rangle\) we denote the standard scalar product
in \(L^2((-a,a))\):%
 \begin{equation*}\textup{For
}u,\,v\in{}L^2((-a,a)),\ \ \ \langle{}u,v\rangle=
\int\limits_{-a}^{a}u(t)\overline{v(t)}\,dt\,.
\end{equation*}
\vspace{1.0ex} \noindent \textsf{The properties of the operators}
\(\mathcal{L}_{{}_\textup{min}}\)
\textsf{and} \(\mathcal{L}_{{}_\textup{max}}\):\\[1.0ex]
\hspace*{2.5ex}1.\
\begin{minipage}[t]{0.90\linewidth}
\textit{The operator \(\mathcal{L}_{{}_\textup{min}}\) is symmetric:
\begin{equation}%
\label{MinSym}%
\langle{}\mathcal{L}_{{}_\textup{min}}x,y\rangle=
\langle{}x,\mathcal{L}_{{}_\textup{min}}y{}\rangle\,,
\quad \forall {}x,\,y\in{}%
\mathcal{D}_{\!{}_{\scriptstyle\mathcal{L}_{\textup{min}}}};
\end{equation}%
In other words, the operator \(\mathcal{L}_{{}_\textup{min}}\) is contained in its adjoint}:
\begin{equation*}%
\mathcal{L}_{{}_\textup{min}}\subseteq(\mathcal{L}_{{}_\textup{min}})^\ast\,;
\end{equation*}%
\end{minipage}{\ }\\[1.0ex]
\hspace*{2.5ex}2.\
\begin{minipage}[t]{0.90\linewidth}
\textit{The operators \(\mathcal{L}_{{}_\textup{min}}\) and \(\mathcal{L}_{{}_\textup{max}}\)
are mutually adjoint}:
\begin{equation}%
\label{ATDO}%
(\mathcal{L}_{{}_\textup{min}})^\ast=\mathcal{L}_{{}_\textup{max}},\quad
(\mathcal{L}_{{}_\textup{max}})^\ast=\mathcal{L}_{{}_\textup{min}}\,;
\end{equation}%
\end{minipage}{\ }\\[1.0ex]

In 1930 John von Neumann, \cite{Neu}, has found a criterion for
the existence of a self-adjoint extension of a symmetric operator
\(A_0\) and has described all such extensions. This criterion is
formulated in terms of deficiency indices of the symmetric
operator.

\begin{definition}
\label{DefDefInd} Let \(A_0\) be an operator in a Hilbert space
\(\mathfrak{H}\). We assume that the domain of definition
\(\mathcal{D}_{A_0}\) is dense in \(\mathfrak{H}\) and that the
operator \(A_0\) is symmetric, that is\,%
\footnote{%
The relation \eqref{SymInc} means that
\(\mathcal{D}_{A_0}\subseteq\mathcal{D}_{A_0^\ast}\) and
\(A_0x=A_0^\ast{}x\,\forall{}x\in\mathcal{D}_{A_0}\).}
\begin{equation}
\label{SymInc}%
 A_0\subseteq{}(A_0)^{\ast}.
\end{equation}

For complex number \(\lambda\), consider the subspace
\begin{subequations}
\begin{equation}
\label{DefSp1}
\mathcal{N}_{\lambda}=\mathfrak{H}\ominus\big((A_0-\lambda{}I)\mathcal{D}_{\!A_0}\big)\,,
\end{equation}
\end{subequations} or, what is equivalent,
\begin{equation}
\label{DefSp}%
 \mathcal{N}_{\lambda}=\lbrace{}x\in\mathfrak{H}:\
 (A_0)^\ast{}x=\overline{\lambda}{}x\rbrace\,.
\end{equation}
 The dimension
\(\dim{}\mathcal{N}_{\lambda}\) is constant in the upper
half-plane \(\textup{Im}\,\lambda>0\) and in the lower half-plane
\(\textup{Im}\,\lambda<0\):
\begin{subequations}
\label{DefInd}
\begin{align}
\label{DefInd+}
\dim{}\mathcal{N}_{\lambda}= n+\,,\quad & \textup{Im}\,\lambda>0,\\
\label{DefInd-} \dim{}\mathcal{N}_{\lambda}= n-\,,\quad &
\textup{Im}\,\lambda<0\,.
\end{align}
\end{subequations}
The numbers \(n^+\) and \(n^{-}\) are said to be the
\emph{deficiency indices of the operator \(A_0\)}, and the
subspace \(\mathcal{N}_{\lambda}\) is said to be \emph{the
deficiency subspace corresponding to the value \(\lambda\)}.
\end{definition}

\begin{nonumtheorem}[von Neumann]{\ } \\ 
\hspace*{1.5ex}\textup{1}.\,%
\begin{minipage}[t]{0.95\linewidth}
The densely defined symmetric operator admits selfadjoint
expansions is and only if its deficiency indices are equal:
\begin{equation}
\label{EqDI}
n_{+}=n_{-}\,.
\end{equation}
\end{minipage}\\
\hspace*{1.5ex}\textup{2}.\,%
\begin{minipage}[t]{0.95\linewidth}
Assume that a symmetric operator \(A_0\) is closed and its
deficiency indices are equal. Choose a pair of non-real conjugated
complex numbers, for example \(\lambda=i\),
\(\overline{\lambda}=-i\). The set of all selfadjoint extensions
of the operator \(A\) is in one-to-one correspondence with the set
of all unitary operators acting from the deficiency subspace
\(\mathcal{N}_i\) into the deficiency subspace
\(\mathcal{N}_{-i}\). In particular, if \(n+=n_-=0\), the operator
\(A_0\) already is selfadjoint.
\end{minipage}
\end{nonumtheorem}
We apply the von Neumann Theorem to the situation where the
operator \(\mathcal{L}_{\textup{min}}\) is taken as the operator
\(A_0\). Then the equation \[(A_0)^\ast{}x=\lambda{}x\] takes the
form
\[\mathcal{L}_{\textup{max}}x=\overline{\lambda}x\,,\] that is the
differential equation
\begin{equation}
\label{DEFEWPr}
-\frac{d\,\,}{dt}\bigg(1-\frac{t^2}{a^2}\bigg)\frac{dx(t)}{dt}+t^2x(t)=
\overline{\lambda}{}x(t),\quad
 t\in(-a,a),
\end{equation}
under the extra condition \(x(t)\in{}L^2(-a,a).\) In particular,
the dimension of the deficiency space \(\mathcal{N}_{\lambda}\)
coincides with the dimension of the linear space of the set of
solutions of the equation \eqref{DEFEWPr} belongings to
\(L^2(-a,a)\). According to Lemma \ref{InArS}, every solution of
the equation \eqref{DEFEWPr} belongs to \(L^2(-a,a)\). Thus we
prove the following
\begin{lemma}
\label{LDefInd}%
 For the operator \(\mathcal{L}_{\textup{min}}\), the deficiency
 indices are:
 \begin{equation}
 \label{dimdo}%
 n_{+}(\mathcal{L}_{\textup{min}})=2,\quad n_{-}(\mathcal{L}_{\textup{min}})=2\,.
\end{equation}
\end{lemma}
Thus, the operator \(\mathcal{L}_{\textup{min}}\) is symmetric,
but not selfadjoint. The set of all its selfadjoint extensions can
by parameterized by the set of all \(2\times{}2\) unitary
operators acting from the two-dimensional Hilbert space
\(\mathcal{N}_{i}\) onto  two-dimensional Hilbert space
\(\mathcal{N}_{-i}\), where \(\mathcal{N}_{\pm{}i}\) are defect
subspaces of the operator \(\mathcal{L}_{\textup{min}}\).
\\[1.0ex]

\noindent
\textsf{Selfadjoint extensions of operators and self-orthogonal subspaces.}\\
J. von Neumann, \cite{Neu},  reduced the construction of a
selfadjoint extension for a symmetric operator \(A_0\) to an
equivalent problem of construction of an unitary extension of an
appropriate isometric operator - the Caley transform of this
symmetric operator. This approach was also developed by M.\,Stone,
\cite{St}, and then used by many others.

In some situations, it is much more convenient to use the
construction of extensions based on the so called \emph{boundary
forms}. Especially convenient is the usage of this construction
for differential operators. The first version of the extension
theory based on abstract symmetric boundary conditions, was
developed by J.W.\,Calkin, \cite{Cal}. Subsequently various
versions of the extension theory of symmetric operators in terms
of abstract boundary conditions were developed. The dual problem
of the descriptions of extensions of symmetric boundary relations
was also considered. See \cite{RoB}, \cite{Koch}, \cite{Br}.

Considering the symmetric operator \(A_0\) \eqref{SymInc} acting
in a Hilbert space \(\mathfrak{H}\), we introduce the bilinear
form form \(\Omega\)
\begin{subequations}
\label{HerFor}
\end{subequations}
\begin{equation}%
\label{HerFor1}
\tag{\ref{HerFor}a}
\Omega(x,y)=
\frac{\langle{}A^{\ast}_0{}x,y\rangle-\langle{}x,A^\ast_0{}y\rangle}{i}\,,
\quad \Omega:\,\mathcal{D}_{\!A^{{}^{\ast}}_{0}}\times{}%
 \mathcal{D}_{\!A^{{}^\ast}_{0}}\to\mathbb{C}\,.
\end{equation}%

The bilinear form \(\Omega\) is hermitian:
\begin{equation*}%
\Omega(x,y)=\overline{\Omega(y,x)},\qquad \forall{}x,y\in \mathcal{D}_{\!A^{{}^{\ast}}_{0}}\,,
\end{equation*}%
and possesses the property
\begin{equation*}%
\Omega(x,y)=0,\quad \forall\,x\in{}\mathcal{D}_{A^{{\ast}}_{0}},\,y\in{}\mathcal{D}_{A_0}\,.
\end{equation*}%
This property allows to consider the form \(\Omega\) as a form on the factor-space \(\mathcal{E}\):
\begin{equation}%
\label{FaSp}%
\mathcal{E}=
\mathcal{D}_{A_0^\ast}\big%
/\mathcal{D}_{A_0}\,.
\end{equation}%
We use the same notation for the form induced on the factor space \(\mathcal{E}\):
\begin{equation}%
\tag{\ref{HerFor}b}
\label{HerFor2}
\Omega(x,y)=\frac{\langle{}A^{\ast}{}x,y\rangle-\langle{}x,A^\ast{}y\rangle}{i}\,,
\quad
\Omega:\,\mathcal{E}\times{}\mathcal{E}\to\mathbb{C}\,.
\end{equation}%
\begin{definition}
\label{BoundForm}
The form \(\Omega\), \eqref{HerFor}, is said to be the \emph{boundary form.}
The factor space \(\mathcal{E}\) is said to be the \emph{boundary space}.
\end{definition}
It turns out that
\begin{equation}%
\label{DimBS}
 \dim{}\mathcal{E}=n_{+}+n_{-}\,,
\end{equation}%
where \(n_{+}\) and \(n_{-}\) are deficiency indices of the operator \(A_0\),
and that \emph{the form \(\Omega\) is not degenerate on \(\mathcal{E}\)}.
The non-degeneracy of the form means:
\begin{equation}%
\textup{ \footnotesize For every non-zero }x\in\mathcal{E} ,
\textup{ \footnotesize  there exists } y\in\mathcal{E}
\textup{ \footnotesize  such that }\Omega(x,y)\not=0\,.
\end{equation}%
Let \(\mathcal{S}\) be a subspace of the factor space \(\mathcal{E}\):
\begin{subequations}
\label{extsp}
\begin{equation}%
\label{extsp1}
\mathcal{S}\subseteq\mathcal{E}\,.
\end{equation}%
 We identify
\(\mathcal{S}\) and its preimage with respect to the factor-mapping
\(\mathcal{D}_{A_0^\ast}\to\mathcal{D}_{A_0^\ast}\big/\mathcal{D}_{A_0}\,(=\mathcal{E})\)
and use  the same notation \(\mathcal{S}\) for a subspace in \(\mathcal{E}\)
and its preimage in \(\mathcal{D}_{A_0}\):
\begin{equation}%
\label{extsp2}
\mathcal{D}_{A_0}\subseteq\mathcal{S}\subseteq\mathcal{D}_{A_0^\ast}\,.
\end{equation}%
\end{subequations}
To every \(\mathcal{S}\) satisfying \eqref{extsp2}, an extension of the
operator \(A_0\) is related. We denote this extension by \(A_{\mathcal{S}}\):
\[\mathcal{D}_{A_{\mathcal{S}}}=\mathcal{S}\,,\quad
A_0\subseteq{}A_{\mathcal{S}}\subseteq{}A_0^\ast\,. \]
The operator \((A_{\mathcal{S}})^{\ast}\), which is the operator adjoint to
the the operator \(A_{\mathcal{S}}\), is related to the subspace
\(\mathcal{S}^{\bot}\):
\begin{equation}%
\label{sadop}
(A_{\mathcal{S}})^\ast=A_{\mathcal{S^\bot}}\,,
\end{equation}%
where \(\mathcal{S}^{\bot_{\Omega}}\) is the \emph{orthogonal
complement} of the subspace \(\mathcal{S}\) with respect to the
hermitiam form \(\Omega\):
\begin{equation}%
\label{orcom}
\mathcal{S}^{\bot_\Omega}=\lbrace{}x\in\mathcal{E}:\ \ \Omega(x,y)=0\ \ %
 \forall\,y\in\mathcal{S}\rbrace\,.
\end{equation}%
In particular we prove the following result:
\begin{lemma}
\label{CrSel}
The extension \(A_{\mathcal{S}}\) of the symmetric operator \(A_0\)
is a selfadjoint operator: \(A_{\mathcal{S}}=(A_{\mathcal{S}})^\ast\),
if and only if the subspace \(\mathcal{S}\) which appears in \eqref{extsp2}
possesses the property:
\begin{equation}%
\label{SeOrSu}
\mathcal{S}=\mathcal{S}^{\bot_\Omega}\,.
\end{equation}%
\end{lemma}
\begin{definition}
\label{DeSoSu} The subspace \(\mathcal{S}\) of the space
\(\mathcal{E}\) is said to be \(\Omega\)-\emph{self-complementary}
if it possess the property \eqref{SeOrSu}.
\end{definition}

It turns out that self-complementary subspaces exist if end only
if the form \(\Omega\), \eqref{HerFor2}, has equal numbers of
positive and negative squares. (Which conditions is equivalent to
the condition \(n_{+}=n_{-}\).)

Thus, \emph{the problem of description of all self-adjoint
extension of a symmetric operator \(A_0\) can be reformulate as
the problem of description of subspaces of the space
\(\mathcal{E}\), \eqref{FaSp}, which are self-complementary with
respect to the (non-degenerated) boundary form \(\Omega\), \eqref{HerFor2}.}\\[1.0ex]

\noindent
\textsf{Selfadjoint extensions of symmetric differential operators.}
The description of selfadjoint extensions of a symmetric operator \(A_0\) becomes
especially transparent in the case when this symmetric operator is
formally selfadjoint ordinary differential operator, regular or singular. In this case the
\emph{hermitian form} \(\Omega\), \eqref{HerFor1}, can be expressed in term of
\emph{boundary conditions} of functions from domain of definitions of
the operator~\(A_0^\ast\). This justifies the terminology introduced in
Definition~\ref{BoundForm}.

We illustrate the situation as applied to the case where the
symmetric operator \(A_0\) is the minimal differential operator
\(\mathcal{L}_{\textup{min}}\) generated by the formal prolate
spheroid differential operator \(L\). Then the adjoint operator
\(A_0^\ast\) is the maximal differential operator
\(\mathcal{L}_{\textup{max}}\) (See Definitions \ref{MinDO} and
\ref{MaxDO}.) The problem of description of selfadjoint
differential operators generated by a given formal differential
operator, has the long history. See, for example, \cite{Kr},
\cite[Chapter 5]{Nai}. The book of \cite{DuSch} is the storage of
wisdom in various aspects of the operator theory, in particular is
self-adjoint ordinary differential operators. See especially
Chapter XIII of \cite{DuSch}.

In principle we may incorporate the question of description of
selfadjoint boundary condition for the prolate spheroid
differential operators in one or other of the existing abstract
schemes which is devoted to such a description in one or other
generality. However to adopt our question to such a scheme one
need to agree the notation, the terminology, etc. This auxiliary
work may obscure the presentation. To make the presentation more
transparent, we prefer to act independently on the existing
general considerations and to develop what we need from the blank
page.

We use the notations
\begin{equation*}%
p(t)=1-\frac{t^2}{a^2},\quad q(t)=t^2, \quad -a<t<a.
\end{equation*}
In this notation, the formal differential operator \(L\) introduced in
\eqref{DiffExp} is:
\begin{equation*}%
(Lx)(t)=-\frac{d}{dt}\bigg(p(t)\frac{dx(t)}{dt}\bigg)+q(t)x(t), \quad -a<t<a\,.
\end{equation*}
For every \(x,y\in\mathcal{A},\)
\begin{equation*}%
(Lx(t))\,\overline{y(t)}-x(t)\,{}(\overline{Ly(t)})=\frac{d}{dt}[x(t),y(t)],\quad -a<t<a\,,
\end{equation*}%
where
\begin{equation*}%
[x(t),y(t)]=-p(t)\bigg(\frac{dx(t}{dt}\overline{y(t)}-x(t)\overline{\frac{dy(t}{dt}}\bigg)\,.
\end{equation*}%
Therefore, for every \(x,y\in\mathcal{A}\) and for every
\(\alpha,\,\beta:\,-a<\alpha<\beta<a\),
\begin{equation}%
\label{LRBFo3}%
\int\limits_{\alpha}^{\beta}\Big((Lx(t))\,%
\overline{y(t)}-x(t)\,{}\big(\overline{Ly(t)}\big)\Big)\,dt=[x,y](\beta)-[x,y](\alpha)\,.
\end{equation}%
\begin{subequations}
\begin{lemma}
\label{EBV}
For every \(x,y\in\mathcal{D}_{\mathcal{L}_{\textup{max}}}\), there exist the
limits
\label{LRBFo}%
\begin{equation}%
\label{LRBFo1}%
[x,y]_{-a}\stackrel{\textup{\tiny def}}{=}\lim_{\alpha\to{}-a+0}[x,y](\alpha),\quad
[x,y]^{a}\stackrel{\textup{\tiny def}}{=}\lim_{\beta\to{}a-0}[x,y](\beta)\,.
\end{equation}%
\end{lemma}
\begin{proof}
Since \(x(t),\,y(t), (Lx(t)), (Ly(t))\) belongs to \(L^2((-a,a))\), then
\(\int\limits_{-a}^{a}\Big|((Lx(t))\,%
\overline{y(t)}-x(t)\,{}\big(\overline{Ly(t)}\big)\Big|\,dt<\infty\).
Therefore
\begin{multline*}
\int\limits_{-a}^{a}\Big(((Lx(t))\,%
\overline{y(t)}-x(t)\,{}\big(\overline{Ly(t)}\big)\Big)\,dt=\\
=\lim_{\substack{\alpha\to{-a+0}\\
\beta\to{\,a-0}}}\int\limits_{\alpha}^{\beta}\Big(((Lx(t))\,%
\overline{y(t)}-x(t)\,{}\big(\overline{Ly(t)}\big)\Big)\,dt.
\end{multline*}
Concerning this and related result see for example \cite[Chapter
10]{HuPy}.
\end{proof}
The boundary form \(\Omega\), constructed from the operator
\(A_0=\mathcal{L}_{\textup{min}}\) according to \eqref{HerFor1},
\begin{equation*}%
\Omega_L(x,y)=\frac{\mathcal{\langle{}L}_{\textup{max}}\,x,y\rangle-
\langle{}x{},\mathcal{L}_{\textup{max}}\,y\rangle}{i}
\end{equation*}%
can be expressed in the term of the generalized boundary values:
\begin{equation}%
\label{LRBFo2}%
\Omega_L(x,y)=\frac{[x,y]^{a}-[x,y]_{-a}}{i}\,\cdot
\end{equation}%
\end{subequations}
According to \eqref{DimBS} and Lemma \ref{LDefInd}, the dimension
of the boundary space \(\mathcal{E}_L\):
\(\mathcal{E}_L=\mathcal{D}_{\mathcal{L}_{\textup{max}}}\big/%
\mathcal{D}_{\mathcal{L}_{\textup{min}}}\) is:
\begin{equation}%
\label{DimFS}
\dim\,\mathcal{E}_L=4\,.
\end{equation}%
To make calculation explicit, we choose a special basis in the
space \(\mathcal{E}_L\) in which the bilinear form \(\Omega_L\) is
reduced to "sum of squares". The asymptotic behavior of solutions
of the equation \(Lx=0\) near the endpoints of the interval
\((-a,a)\), described in Lemma \ref{ABSNS}, prompts us the choice
of such a basis. Let \(\varphi_{-}(t), \psi_{-}(t),
\varphi_{+}(t), \psi_{+}(t)\) be smooth functions such that
\begin{alignat}{4}
&\varphi_{-}(t)=1,\ &-a<t<-a/2,&\quad &\varphi_{-}(t)&=0, \
&a/2<t<a\,,\notag{}\\
&\psi_{-}(t)=\ln(a+t),\ &-a<t<-a/2,&\quad &\psi_{-}(t)&=0, \
&a/2<t<a\,,\notag{}\\
 &\varphi_{+}(t)=0, \ &-a<t<-a/2\,,&\quad &\varphi_{+}(t)&=1,\
&a/2<t<a\,,\notag{}\\
&\psi_{+}(t)=0, \  &-a<t<-a/2\,,&\quad &\psi_{+}(t)&=\ln(a-t),\
&a/2<t<a\,.
\label{BounF}
\end{alignat}

It is cleat that if \(\chi\) is an arbitrary smooth
real valued function, then \(\Omega_L(\chi,\chi)=0)\). In
particular,
\begin{subequations}
\label{CBF}
\begin{equation}
\label{CBF1}
\Omega_L(\chi,\chi)=0,\ \ \textup{if \(\chi\) is one of the functions } \varphi_{-},\psi_{-},
\varphi_{+},\psi_{+}\,.
\end{equation}
It is clear that
\begin{equation}
\label{CBF2}
\Omega_L(\chi_{-},\chi_{+})=0,\ \ \textup{if \(\chi_{\pm}\) is one of the functions }
\varphi_{\pm},\psi_{\pm}\,.
\end{equation}
 Direct calculation shows that
\begin{equation}
\label{CBF3}
\Omega_L(\varphi_{-},\psi_{-})=-\frac{2}{a},\quad \Omega_L(\varphi_{+},\psi_{+})=-\frac{2}{a}\,.
\end{equation}
\end{subequations}
Thus, the Gram matrix (with respect to the hermitian form \(\Omega_L\))
 of the vectors \(\varphi_{-}\), \(\psi_{-}\),
 \(\varphi_{+}\), \(\psi_{+}\) is:
 \begin{equation}
 \label{GrMa}
 \frac{a}{2}  \cdot%
 \begin{bmatrix}
\Omega_L(\varphi_{-},\varphi_{-})&\Omega_L(\varphi_{-},\psi_{-})%
&\Omega_L(\varphi_{-},\varphi_{+})&\Omega_L(\varphi_{-},\psi_{+})\\
\Omega_L(\psi_{-},\varphi_{-})&\Omega_L(\psi_{-},\psi_{-})%
&\Omega_L(\psi_{-},\varphi_{+})&\Omega_L(\psi_{-},\psi_{+})\\
\Omega_L(\varphi_{+},\varphi_{-})&\Omega_L(\varphi_{+},\psi_{-})%
&\Omega_L(\varphi_{+},\varphi_{+})&\Omega_L(\varphi_{+},\psi_{+})\\
\Omega_L(\psi_{+},\varphi_{-})&\Omega_L(\psi_{+},\psi_{-})%
&\Omega_L(\psi_{+},\varphi_{+})&\Omega_L(\psi_{+},\psi_{+})
 \end{bmatrix}=J,
 \end{equation}
 where
 \begin{equation}
 \label{IndMe}
 J=
\begin{bmatrix}
0 &i & 0&0 \\
 -i & 0& 0&0 \\
 0 & 0& 0& i\\
 0& 0& -i&0
\end{bmatrix}\,.
\end{equation}
The rank of the Gram matrix is is equal to the dimension of the space~\(\mathcal{E}_L\):
\[\textup{rank}\,J=\dim{}\mathcal{E}_L=4\,.\]
Therefore, the vectors \(\varphi_{-}\), \(\psi_{-}\),
 \(\varphi_{+}\), \(\psi_{+}\) generate the space
 \(\mathcal{E}_L=\mathcal{D}_{\mathcal{L}_{\textup{max}}}\big/%
\mathcal{D}_{\mathcal{L}_{\textup{min}}}\). In particular,
the domain of definition \(\mathcal{D}_{\mathcal{L}_{\textup{min}}}\) of the minimal
differential operator \(\mathcal{L}_{\textup{min}}\) can be characterized by means
of the boundary conditions:
\begin{multline}
\label{DoDeMiO}
\mathcal{D}_{\mathcal{L}_{\textup{min}}}=\big\lbrace{}x(t):
x(t)\in\mathcal{D}_{\mathcal{L}_{\textup{max}}},\ \textup{and}\\
\,[x,\varphi_{-}]_{-a}=0,\,[x,\psi_{-}]_{-a}=0,\,
[x,\varphi_{+}]^{a}=0,\,[x,\psi_{+}]^{a}=0\big\rbrace\,,
\end{multline}
where the forms \([\,\,,\,\,]_{-a},\,[\,\,,\,\,]^a\) were
introduced in \eqref{LRBFo}.
\begin{lemma}
\label{OJOrt}%
Let \(\Omega_L\) be a bilinear form in the space \(\mathcal{E}\) defined by
\eqref{LRBFo}, and \(J\) be the matrix \eqref{IndMe}.

The vector \(x^{1}=\alpha_{-}^{1}\varphi_{-}+\beta_{-}^{1}\psi_{-}+
\alpha_{+}^{1}\varphi_{+}+\beta_{+}^{1}\psi_{+}\in\mathcal{E}_L\) is \(\Omega_L\)~-orthogonal
to the vector \(x^{2}=\alpha_{-}^{2}\varphi_{-}+\beta_{-}^{2}\psi_{-}+
\alpha_{+}^{2}\varphi_{+}+\beta_{+}^{2}\psi_{+}\in\mathcal{E}_L\), that is
\begin{subequations}
\label{Ort}
\begin{equation}%
\label{Ort1}
\Omega_L(x^1,x^2)=0,
\end{equation}
 if and only
if the vector-row \(v_{x^1}=[\alpha_{-}^{1},\beta_{-}^{1},\alpha_{+}^{1},\beta_{+}^{1}]\in\mathcal{V}\)
is \(J\)-orthogonal to the vector-row \(v_{x^2}=[\alpha_{-}^{2},\beta_{-}^{2},\alpha_{+}^{2},\beta_{+}^{2}]\in\mathcal{V}\), that is
\begin{equation}%
\label{Ort2}
v_{x^1}J\,v_{x^2}^{\,\ast}=0\,,
\end{equation}
\end{subequations}
where \(\mathcal{V}\) is the space \(\mathbb{C}^4\) of vector-rows equipped by the standard hermitian metric, and the star \(\ast\) is the Hermitian conjugation.
\end{lemma}
Thus, the problem of description of self-complementary extensions
of the operator \(\mathcal{L}_{\textup{min}}\) is equivalent to
the problem of description
of \(\Omega_L\)-self-complementary\,%
\footnote{%
As soon as the notion of \(J\)-orthogonality of two vectors is introduced,
\eqref{Ort2}, the notions of \(J\)-orthogonal complement and \(J\)-self-orthogonal
subspaces can be introduced as well.
} %
subspaces in \(\mathcal{E}\), which in its turn is equivalent to
the problem of description of \(J\)-self-complementary subspaces
in \(\mathbb{C}^4\). The last problem is a problem of the
indefinite linear algebra and admits an explicit  solutions. We
set
\begin{subequations}%
\label{Jpro}
\begin{equation}
\label{Jpro1}
P_{+}=\frac{1}{2}(I+J),\quad P_{-}=\frac{1}{2}(I-J)\,,
\end{equation}
More explicitly,
\begin{equation}
\label{Jpro2}
P_{+}=\frac{1}{2}
\begin{bmatrix}
1 &i & 0&0 \\
 -i & 1& 0&0 \\
 0 & 0& 1& i\\
 0& 0& -i&1
\end{bmatrix}\,,\quad
P_{-}=\frac{1}{2}
\begin{bmatrix}
1 &-i & 0&0 \\
 i & 1& 0&0 \\
 0 & 0& 1& -i\\
 0& 0& i&1
\end{bmatrix}\,.
\end{equation}
\end{subequations}
The matrix \(J\), \eqref{IndMe}, possesses the properties
\begin{equation*}
J=J^{\ast},\quad J^2=I.
\end{equation*}
Therefore the matrices \(P_{+},\ P_{-}\), \eqref{Jpro1},  possess  the properties
\begin{align}
\label{PrPr}
P_{+}^2=P_{+},\quad P_{-}^2&=P_{-}, \quad P_{+}=P_{+}^{\,\ast},\quad
\quad P_{-}=P_{-}^{\,\ast},\\
  P_{+}P_{-}&=0, \quad P_{+}+P_{-}=I\,.
\end{align}
In other words, the matrices \(P_{+},\ P_{-}\) are orthogonal projector matrices.
These matrices project the space \(\mathcal{V}\) onto subspaces
\(\mathcal{V}_{+}\) and \(\mathcal{V}_{-}\):
\begin{equation}
\label{PMSu}%
\mathcal{V}_{+}=\mathcal{V}P_{+},\ \mathcal{V}_{-}=\mathcal{V}P_{-}\,.
\end{equation}
These  subspaces are orthogonally complementary:
\begin{equation}
\label{PMSu1}%
\mathcal{V}_{+}\oplus\mathcal{V}_{-}=\mathcal{V}\,.
\end{equation}
The vector rows
\begin{subequations}
\label{bas}
\begin{alignat}{2}
\label{bas1}
e^1_{+}&=[1,\phantom{-}i,0,0],\quad& e^2_{+}&=[0,0,1,\phantom{-}i]\\
\intertext{and}
\label{bas2}
e^1_{-}&=[1,-i,0,0],\quad & e^2_{-}&=[0,0,1,-i]%
\end{alignat}
\end{subequations}
form orthogonal\,%
\footnote{\,%
 In the standard scalar product on
\(\mathcal{V}=\mathbb{C}^4\).} %
 bases in \(\mathcal{V}_{+}\) and \(\mathcal{V}_{-}\) respectively.

 It turns out that \(J\)-self-orthogonal subspaces of the space \(\mathcal{V}\)
are in one-to-one correspondence with unitary operators acting from
\(\mathcal{V}_{+}\) onto \(\mathcal{V}_{-}\).
\begin{definition}
Let \(U\) be an unitary operator  acting from \(\mathcal{V}_{+}\) onto
\(\mathcal{V}_{-}\). As the vector-row \(v\) runs  over the whole subspace
\(\mathcal{V}_{+}\), the vector \(v+vU\) runs over a subspace
of the space \(\mathcal{V}\). This subspace is denoted by \(\mathcal{S}_U\):
\begin{equation}
\mathcal{S}_U=\big\lbrace{}v+vU\big\rbrace,\quad\textup{where }v%
\textup{ runs over the whole }\mathcal{V}_{+}\,.
\end{equation}
\end{definition}
\begin{lemma}{\ }\\
\label{DJSO} \hspace*{1.5ex}\textup{1.}
\begin{minipage}[t]{0.95\linewidth} Let \(U\) be an unitary
operator  acting from \(\mathcal{V}_{+}\) onto
\(\mathcal{V}_{-}\). Then the subspace \(\mathcal{S}_U\)  is
\(J\)-self-complementary, that is
\[\mathcal{S}_U=\mathcal{S}_U^{\bot_J}\,.\]
\end{minipage}\\[1.0ex]
\hspace*{1.5ex}\textup{2.} \begin{minipage}[t]{0.95\linewidth}
 Every \(J\)-self-complementary subspace \(\mathcal{S}\)
of the space \(\mathcal{V}\) is of the form \(\mathcal{S}_U\):
\[\mathcal{S}=\mathcal{S}_U\]
for some unitary operator \(U:\,\mathcal{V}_{+}\to\mathcal{V}_{+-}\).
\end{minipage}\\[1.0ex]
\hspace*{1.5ex}\textup{3.} \begin{minipage}[t]{0.95\linewidth}
 The correspondence between \(J\)-self-complementary subspaces and
unitary operators acting from \(\mathcal{V}_{+}\) onto
\(\mathcal{V}_{-}\) is one-to-one;
\[(U_1=U_2)\Leftrightarrow{}(\mathcal{S}_{U_1}=\mathcal{S}_{U_2})\,.\]
\end{minipage}
\end{lemma}
\begin{proof}
\textbf{1}.\,The mapping \(v\to{}v+Uv\) is one-to-one mapping from \(\mathcal{V}_{+}\)
onto \(\mathcal{S}_U\). Indeed, this mapping is surjective by definition of
the subspace \(\mathcal{S}_U\). This mapping is also injective.
 The equality \(v+Uv=0\) implies that \(v=Uv=0\) since\,%
 \footnote{%
 Recall that
 \(v\in\mathcal{V}_{+},\,Uv\in\mathcal{V}_{-}\), %
and \(\mathcal{V}_{+}\bot \mathcal{V}_{-}\).
 } %
  \(v\bot\,Uv\). In particular, \(\dim\mathcal{S}_U=\dim\mathcal{V}_{+}\,(=2)\).

If \(v_1\) and \(v_2\) are two arbitrary vectors from \(\mathcal{V}_{+}\),
then the vectors \(w_1=v_1+v_1U\) and \(w_2=v_2+v_2U\) are \(J\)-orthogonal:
\(w_1Jw_2^{\ast}=0\). Indeed, since \(J=P_{+}-P_{-}\) and %
\(v_k=v_kP_{+}, v_kU=v_kUP_{-}\,,k=1,2\), then, using the properties \eqref{PrPr}
of \(P_{+}\) and  \(P_{-}\), we oobtain
\begin{align*}%
w_1Jw_2^{\ast}=(v^1P_{+}&+v^1UP_{-})(P_{+}-P_{-})
(P_{+}^\ast{}v_2^\ast+P_{-}^\ast{}U^\ast{}v_2^\ast)=\\
&=v_1v_2^\ast{}-v_1UU^\ast{}v_2^\ast\,.
\end{align*}
Since the unitary operator \(U\) preserves the scalar product, then
\(v_1v_2^{\ast}=v_1UU^\ast{}v_2^{\ast}\), hence \(w_1Jw_2^{\ast}=0\).
Thus, \(\mathcal{S}_U\subseteq(\mathcal{S}_U)^{\bot_J}\).
(The symbol \(\bot_J\) means \(J\)-orthogonal complement.)
Since the Hermitian form \((v_1,v_2)\to{}v_1Jv_2^\ast\) is non-degenerate
on \(\mathcal{V}\), then \(\dim(\mathcal{S}_U^{\bot_J})=
\dim{\mathcal{V}}-\dim{\mathcal{S}_U}\). Because %
\(\dim{\mathcal{V}}-\dim{\mathcal{S}_U}=\dim{\mathcal{S}_U}\), we
have \(\dim{\mathcal{S}_U}=\dim(\mathcal{S}_U^{\bot_J})\). Hence,
\(\mathcal{S}_U=(\mathcal{S}_U)^{\bot_J}\), i.e. the subspace
\(\mathcal{S}_U\) is \(J\)-self-complementary.

\textbf{2}.\,Let \(\mathcal{S}\) be a \(J\)-self-orthogonal subspace.
If \[v\in\mathcal{S},\,v=v_{+}+v_{-},\,v_{\pm}\in\mathcal{V}_{\pm}\,,\]
then the condition \(v{\bot_J}v=0\), that is the condition \(vJv^\ast=0\)
means that \(v_1v_1^\ast=v_2v_2^\ast\). Therefore, if \(v_1=0\), then also
\(v=0\). This means that the projection mapping \(v\to{}vP_{+}\), considered
as a mapping from \(\mathcal{S}\to\mathcal{V}_{+}\), is injective.
For \(J\)-self-orthogonal subspace \(\mathcal{S}\) of the space \(\mathcal{V}\),
the equality \(\dim\mathcal{S}=\dim\mathcal{V}-\dim\mathcal{S}\) holds.
Hence \(\dim\mathcal{S}=\dim\mathcal{V}_{+}\). Therefore, the injective linear
mapping \(v\to{}P_{+}\) is surjective. The inverse mapping
is defined on the whole subspace \(\mathcal{V}_{+}\) and can by presented
in the form \(v=v_1+v_1U\), where \(U\) is a linear operator acting
from \(\mathcal{V}_{+}\) into \(\mathcal{V}_{-}\). This  mapping \(v_1\to{}v_1+v_1U\)
maps the subspace \(\mathcal{V}_{+}\) onto the subspace \(\mathcal{S}\).

 Since \(vJv^\ast=0\),
then \(v_1v_1^\ast=v_2v_2^\ast\), where \(v_2=v_1U\). Since
\(v_1\in\mathcal{V}_{+}\) is arbitrary, this means that the
operator \(U\) is isometric. Since
\(\dim{}\mathcal{V}_{+}=\dim{}\mathcal{V}_{-}\), the operator
\(U\) is unitary. Thus, the originally given
\(J\)-self-complementary subspace \(\mathcal{S}\) is of the form
\(\mathcal{S}_U\), where \(U\) is an unitary operator acting from
\(\mathcal{V}_{+}\) to \(\mathcal{V}_{-}\).

The coincidence \(\mathcal{S}_{U_1}=\mathcal{S}_{U_2}\) means that
every vector of the form \(v_1+v_1U_1\), where \(v_1\in\mathcal{V}_{+}\)
can also be presented in the form \(v_2+v_2U_12\) with some \(v_2\in\mathcal{V}_{+}\):
\[v_1+v_1U_1=v_2+v_2U_2\,.\]
Since \(v_1,\,v_2\in\mathcal{V}_{+},\ v_1U_1,\,v_1U_2\in\mathcal{V}_{-}\),
then \(v_1=v_2\), and \(v_1U_1=v_1U_2\). The equality \(v_1U_1=v_1U_2\)
for every \(v_1\in\mathcal{V}_{+}\) means that \(U_1=U_2\).
Thus, \((\mathcal{S}_{U_1}=\mathcal{S}_{U_2})\Rightarrow(U_1=U_2)\).
\end{proof}
Choosing the orthogonal bases \eqref{bas} in the subspaces \(\mathcal{V}_{+}\)
and \(\mathcal{V}_{+}\), we represent an unitary operator \(U\) by the appropriate
unitary matrix:
\begin{alignat*}{2}
e^1_{+}U&=&\ e^1_{-}u_{11}&+e^2_{-}u_{21},\\
e^2_{+}U&=&\ e^1_{-}u_{12}&+e^2_{-}u_{22}.
\end{alignat*}
The
following result is a  reformulation of Lemma \ref{DJSO}:
\begin{lemma}
\label{BSOS}%
Let \(\mathcal{V}\) be the space \(\mathbb{C}^4\) of \(four\) vector-rows,
\(J\) be a matrix of the form \eqref{IndMe}. With every \(2\times2\) matrix %
\( U=\|u_{pq}\|_{1\leq{}p,q\leq{}2}\), we
associate the pair of vectors \(v^1(U),\,v^2(U)\):
\begin{subequations}
\label{BSOSu}
\begin{alignat}{2}
\label{BSOSu1}
v^1(U)&=e^1_{+}+&\ e^1_{-}u_{11}&+e^2_{-}u_{21},\\
\label{BSOSu2}
v^2(U)&=e^2_{+}+&\ e^1_{-}u_{12}&+e^2_{-}u_{22},
\end{alignat}
\end{subequations}
where \(e^{\,k}_{\pm},\,k=1,2,\)
are the vector-rows of the form \eqref{bas},
and the subspace \(\mathcal{S}_U\) of \(\mathcal{V}\) which is the linear
hull of the vectors \(v^1(U),\,v^2(U)\),
\begin{equation*}
\mathcal{S}_U=\textup{hull}(v^1(U),\,v^2(U))\,.
\end{equation*}
\hspace*{1.5ex}\textup{1}.\,
\begin{minipage}[t]{0.90\linewidth}
If the matrix \(U\) is unitary, then the vectors
\(v^1(U),\,v^2(U)\) are linearly independent, and the subspace
\(\mathcal{S}_U\) is \(J\)-self-complementary.
\end{minipage}\\[1.1ex]
\hspace*{1.5ex}\textup{2}.\,
\begin{minipage}[t]{0.90\linewidth}
Let \(\mathcal{S}\) be a \(J\)-self-complementary subspace of the
space \(\mathcal{V}\). Then \(\mathcal{S}=\mathcal{S}_U\) for some
an unitary matrix \(U\).
\end{minipage}\\[1.1ex]
\hspace*{1.5ex}\textup{3}.\,
\begin{minipage}[t]{0.90\linewidth}
For  unitary matrices \(U_1,\,U_2\),
\begin{equation*}
(\mathcal{S}_{U_1}=\mathcal{S}_{U_2})\Leftrightarrow(U_1=U_2)\,.
\end{equation*}
\end{minipage}
\end{lemma}
The "coordinate" form of the  vectors \(v^1(U),\,v^2(U)\) is:
\begin{alignat}{5}
\label{Coorv}
v^1(U)=%
&\big[&1+u_{11}&,&\,\,\,i(1-u_{11})&,&\,u_{21}\,\,\,\,&,&-iu_{21}\,\,\,\, &\big],
\notag\\
\raisetag{100pt}
v^2(U)=%
&\big[&u_{12}\,\,\,\,&,&\,-iu_{12}\,\,\,\,&,&\,\,1+u_{22}&,&\,\,i(1-u_{22})&\big].
\raisetag{100pt}
\end{alignat}
Remembering, see Lemma \ref{OJOrt}, how \(J\)-self-complementary
subspaces of the space \(\mathcal{V}\) are related to
\(\Omega_L\)-self-complementary subspaces the space
\(\mathcal{E}_L=
\mathcal{D}_{\mathcal{L}_{\textup{max}}}\big/%
\mathcal{D}_{\mathcal{L}_{\textup{min}}}\)
we formulate the following result
\begin{lemma}%
\label{ESoS}%
Let us associate with every \(2\times2\) matrix \( U=\|u_{pq}\|_{1\leq{}p,q\leq{}2}\)
 the pair of vectors \(d^1(U),\,d^2(U)\) of the space \(\mathcal{E}_L\):
 \begin{subequations}%
 \label{VGeE}%
 \begin{align}%
 \label{VGeE1}%
 d^1(U)=&\,(1+u_{11})\varphi_{-}+i(1-u_{11})\psi_{-}+u_{21}\varphi_{+}-iu_{21}\psi_{+}\,,\\
 \label{VGeE2}%
 d^2(U)=&\,u_{12}\varphi_{-}-iu_{12}\psi_{-}+(1+u_{22})\varphi_{+}+i(1-u_{22}\psi_{+}\,,
 \end{align}%
 \end{subequations}%
 where the functions \(\varphi_{\pm},\,\psi_{\pm}\) are defined in \eqref{BounF}.
The subspace \(\mathcal{G}_U\) of the space \(\mathcal{E}_L\) is defined as
the linear hull of the vectors \(d^1(U),\,d^2(U)\):
\begin{equation}
\label{SSU}
\mathcal{G}_U=\textup{hull}\,(d^1(U),\,d^2(U))\,.
\end{equation}
\hspace*{1.5ex}\textup{1}.\,
\begin{minipage}[t]{0.90\linewidth}
If the matrix \(U\) is unitary, then  the subspace
\(\mathcal{S}=\mathcal{G}_U\) is \(\Omega_L\)-self-complementary.
\end{minipage}\\[1.1ex]
\hspace*{1.5ex}\textup{2}.\,
\begin{minipage}[t]{0.90\linewidth}
Let \(\mathcal{S}\) be a \(\Omega_L\)-self-complementary subspace
of the space \(\mathcal{E}_L\). Then \(\mathcal{S}=\mathcal{G}_U\)
for some  an unitary matrix \(U\).
\end{minipage}\\[1.1ex]
\hspace*{1.5ex}\textup{3}.\,
\begin{minipage}[t]{0.90\linewidth}
For  unitary matrices \(U_1,\,U_2\),
\begin{equation*}
(\mathcal{G}_{U_1}=\mathcal{G}_{U_2})\Leftrightarrow(U_1=U_2)\,.
\end{equation*}
\end{minipage}
\end{lemma}%
It is clear that \emph{a subspace
\(\mathcal{S}\subseteq\mathcal{E}_L\) is an
\(\Omega_L\)-self-complementary subspace if and only if its
\(\Omega_L\)-orthogonal complement
\(\mathcal{S}^{\bot_{\Omega_L}}\) is an
\(\Omega_L\)-self-complementary subspace.} The subspace
\((\mathcal{S}_U)^{\bot_{\Omega_L}}\) can be described as:
\[(\mathcal{S}_U)^{\bot_{\Omega_L}}=\big\lbrace{}x\in\mathcal{E}_L:
\Omega_{L}(x,d^1(U))=0,\,\Omega_{L}(x,d^2(U))=0\big\rbrace\,,\]
where \(d^1,\,d^2\) are defined in \eqref{VGeE}, \eqref{BounF}.
Thus Lemma \ref{ESoS} can be reformulated in the following way:%
\begin{lemma}
\label{ODSS}
Let us associate
the pair of vectors \(d^1(U),\,d^2(U)\)
with every \(2\times2\) matrix \( U=\|u_{pq}\|_{1\leq{}p,q\leq{}2}\)
by \eqref{VGeE}, \eqref{BounF}.
The subspace  \(\mathcal{O}_U\) is defined as
\begin{equation}%
\label{OODop}%
\mathcal{O}_U=\big\lbrace{}x\in\mathcal{E}_L:
\,\Omega_{L}(x,d^1(U))=0,\,\Omega_{L}(x,d^2(U))=0\big\rbrace\,.
\end{equation}%
\hspace*{1.5ex}\textup{1}.\,
\begin{minipage}[t]{0.90\linewidth}
If the matrix \(U\) is unitary, then  the subspace
\(S=\mathcal{O}_U\) is \(\Omega_L\)-self-complementary.
\end{minipage}\\[1.1ex]
\hspace*{1.5ex}\textup{2}.\,
\begin{minipage}[t]{0.90\linewidth}
Let \(\mathcal{S}\) be a \(\Omega_L\)-self-complementary subspace
of the space \(\mathcal{E}_L\). Then \(\mathcal{S}=\mathcal{O}_U\)
for some  an unitary matrix \(U\).
\end{minipage}\\[1.1ex]
\hspace*{1.5ex}\textup{3}.\,
\begin{minipage}[t]{0.90\linewidth}
For  unitary matrices \(U_1,\,U_2\),
\begin{equation*}
(\mathcal{O}_{U_1}=\mathcal{O}_{U_2})\Leftrightarrow(U_1=U_2)\,.
\end{equation*}
\end{minipage}
 \end{lemma}

 Thus there is one-to-one correspondence between the set of all \(2\times2\) unitary
 matrices \(U=\|u_{pq}\|_{1\leq{}p,q\leq{}2}\) and the set of all
 \(\Omega_L\)-self-complementary
 subspaces \(\mathcal{S}\) of the space \(\mathcal{E}_L=
\mathcal{D}_{\mathcal{L}_{\textup{max}}}\big/%
\mathcal{D}_{\mathcal{L}_{\textup{min}}}\). This correspondence is described as
\begin{equation}%
\label{Corre}%
\mathcal{S}=\mathcal{O}_U,
\end{equation}%
where \(\mathcal{O}_U\) is defined in \eqref{OODop}, \eqref{VGeE}, \eqref{BounF}.

On the other hand, the subspaces of the space  \(\mathcal{E}_L=
\mathcal{D}_{\mathcal{L}_{\textup{max}}}\big/%
\mathcal{D}_{\mathcal{L}_{\textup{min}}}\) which are
self-complementary with respect to the Hermitian form
\(\Omega_L\), \eqref{LRBFo}, are in one-to-one correspondence to
self-adjoint differential operators generated by the formal
differential operator \(L\), \eqref{DiffExp}. Every self-adjoint
differential operators \(\mathcal{L}\) generated by the formal
differential operator \(L\) is the \emph{restriction} of the
maximal differential operator   \(\mathcal{L}_{\textup{max}}\),
\eqref{MaxDO}, on the appropriate domain of definition. According
to Lemma~\ref{CrSel}, as applied to the operators
\(A_0=\mathcal{L}_{\textup{min}},\,A_0^{\ast}=\mathcal{L}_{\textup{max}}\),
the domains of definition of a selfadjoint extension
\(\mathcal{L}\) of the operator \(\mathcal{L}_{\textup{min}}\) are
those subspaces \(\mathcal{S}\):
\begin{equation}
\mathcal{D}_{\mathcal{L}_{\textup{min}}}\subseteq\mathcal{S}%
\subseteq\mathcal{D}_{\mathcal{L}_{\textup{max}}}
\end{equation}
which are self-complementary with respect to the Hermitian form
\(\Omega_L\), \eqref{LRBFo}. According to Lemma \ref{ODSS},
\(\Omega_L\)-self-complementary subspaces \(\mathcal{S}\) can be
described by means of the conditions
\begin{equation}
\mathcal{S}=\big\lbrace{}x(t)\in\mathcal{D}_{\mathcal{L}_{\textup{max}}}:\,\,%
\Omega_L(x,d^1(U))=0,\,\,\Omega_L(x,d^2(U))=0\big\rbrace\,,
\end{equation}
where \(d^1(U),\,d^2(U)\) are the same that in \eqref{VGeE},\,\eqref{BounF},
\(U\) is an unitary \(2\times2\) matrix.

The conditions \(\Omega_L(x,d^1(U))=0,\,\,\Omega_L(x,d^2(U))=0\)
may be interpreted as a \emph{boundary conditions} posed on functions
\(x\in\mathcal{D}_{\mathcal{L}_{\textup{max}}}\). Let us present these conditions
in more traditional form.

Let as introduce the following notations:
\begin{gather}%
b_{-a}(x)=\lim_{t\to{}{-a+0}}(t+a)\frac{dx(t)}{dt},\quad
b_{a}(x)=\lim_{t\to{}{a-0}}(t-a)\frac{dx(t)}{dt},\notag\\
c_{-a}(x)=\lim_{t\to{}{-a+0}}\bigg((t+a)\ln(a+t)\frac{dx(t)}{dt}-x(t)\bigg),%
\label{GBCo}\\
c_{a}(x)=\lim_{t\to{}{a-0}}\bigg((t-a)\ln(a-t)\frac{dx(t)}{dt}-x(t)\bigg)\,.\notag
\end{gather}
\begin{remark}%
\label{GeBV}
 The values \(b_{-a}(x),\,c_{-a}(x)\) and
\(b_{a}(x),\,c_{a}(x)\) may be considered as generalized boundary
values related to the function
\(x(t)\in\mathcal{D}_{\mathcal{L}_{\textup{max}}}\) at the
endpoints \(-a\) and \(a\) of the interval \((-a,a)\).
\end{remark}
\begin{remark}%
\label{BVBas}
The solutions \(x_{1}^{-},\,x_{2}^{-},\,(x_{1}^{+},\,x_{2}^{+}\) of the
equation \(Lx=\lambda{}x\), which
appears in Lemma \ref{ABSNS}, satisfy the conditions
\begin{alignat*}{4}%
b_{-a}(x_{1}^{-})&=0,\  &c_{-a}(x_{1}^{-})&=-1;\qquad%
&b_{-a}(x_{2}^{-})&=1,\  {}&c_{-a}(x_{2}^{-})&=0\,,\\
b_{a}(x_{1}^{+})&=0,\  &c_{a}(x_{1}^{+})&=-1;\qquad%
&b_{a}(x_{2}^{+})&=1,\  &c_{a}(x_{2}^{+})&=0\,.
\end{alignat*}
\end{remark}
\begin{lemma}%
\label{GBCL}%
For \(x(t)\in\mathcal{D}_{\mathcal{L}_{\textup{max}}}\),
 the limits \eqref{GBCo} exist, are finite, and
\begin{subequations}
\label{GBCTr}
\begin{alignat}{2}%
\label{GBCTr1}
b_{-a}(x)&=\frac{ia}{2}\,\Omega_L(x,\varphi_{-}),&\quad
c_{-a}(x)&=\frac{ia}{2}\,\Omega_L(x,\psi_{-}),\\[1.0ex]
\label{GBCTr2}
b_{a}(x)&=\frac{ia}{2}\,\Omega_L(x,\varphi_{+}),&\quad
c_{a}(x)&=\frac{ia}{2}\,\Omega_L(x,\psi_{+}),
\end{alignat}%
\end{subequations}
where the functions \(\varphi_{\pm}, \psi_{\pm}\) are defined in \eqref{BounF},
and the form \(\Omega_L\) is defined by \eqref{LRBFo}.
\end{lemma}%
\begin{proof} The existence of the limits in \eqref{GBCTr} follows from
Lemma \ref{EBV} applied to the functions \(x(t)\) and \(y(t)=\varphi_{\pm}(t)\)
or \(y(t)=\psi_{\pm}(t)\). The equalities \eqref{GBCTr} can be obtained by
the direct computation using the explicit expressions \eqref{BounF} for the
functions \(\varphi_{\pm}(t),\,\psi_{\pm}(t)\).
\end{proof}
Due to \eqref{GBCTr}, the equality \eqref{GrMa} can be rewritten
as
\begin{equation}
 \label{Rewr}
 \begin{bmatrix}
b_{-a}(\varphi_{-})&c_{-a}(\varphi_{-})&b_{-a}(\varphi_{-})&c_{a}(\varphi_{-})\\
b_{-a}(\psi_{-})&c_{-a}(\psi_{-})&b_{-a}(\psi_{-})&c_{a}(\psi_{-})\\
b_{-a}(\varphi_{+})&c_{-a}(\varphi_{+})&b_{-a}(\varphi_{+})&c_{a}(\varphi_{+})\\
b_{-a}(\psi_{+})&c_{-a}(\psi_{+})&b_{-a}(\psi_{+})&c_{a}(\psi_{+})
 \end{bmatrix}=
\begin{bmatrix}
0&-1&\phantom{-}0&\phantom{-}0\\
1&\phantom{-}0&\phantom{-}0&\phantom{-}0\\
0&\phantom{-}0&\phantom{-}0&-1\\
0&\phantom{-}0&\phantom{-}1&\phantom{-}0
 \end{bmatrix}\,.
 \vspace*{0.8ex}
 \end{equation}
 \begin{remark}
 \label{ChDD}
The characterization \eqref{DoDeMiO} of the domain of definition
\(\mathcal{D}_{\mathcal{L}_{\textup{min}}}\) of the minimal
operator \(\mathcal{L}_{\textup{min}}\) can be presented as:
\begin{multline}
\label{DoDeMiOr}
\mathcal{D}_{\mathcal{L}_{\textup{min}}}=\big\lbrace{}x(t):
x(t)\in\mathcal{D}_{\mathcal{L}_{\textup{max}}},\ \textup{and}\\
\,b_{-a}(x)=0,\,c_{-a}(x)=0,\,
b_{a}(x)=0,\,c_{a}(x)=0\big\rbrace\,,
\end{multline}
\end{remark}%

\vspace{2.0ex}
 In view
of \eqref{GBCTr}, the equalities
\(\Omega_L(x,d^1(U)=0,\,\Omega_L(x,d^2(U)=0)\) take the form
\begin{subequations}
\label{SBoCo}
\begin{align}
\label{SBoCo1}
(1+u_{11})\,b_{-a}(x)-i(1-u_{11})\,c_{-a}(x)+u_{12}\,b_{a}(x)+iu_{12}\,c_{a}(x)&=0\,,\\
\label{SBoCo2}
u_{21}\,b_{-a}(x)+iu_{21}\,c_{-a}(x)+(1+u_{22})\,b_{a}(x)-i(1-u_{22})\,c_{a}(x)&=0
\end{align}
\end{subequations}
\begin{remark}
\label{CUM}%
Since the form \(\Omega_L(x,y)\) is \emph{anti}linear with respect
to the argument \(y\): \(\Omega_L(x,\mu{}y)=
\overline{\mu}\,\Omega_L(x,y)\) for \(\mu\in\mathbb{C}\), the
numbers \(i,-i\) which occurs in \eqref{VGeE} must be replaced
with the numbers \(-i,i\) in appropriate positions in the equality
\eqref{SBoCo}. For the same reason, the numbers \(u_{pq}\) which
occurs in \eqref{VGeE} must be replaced with the numbers
\(\overline{u_{pq}}\) in \eqref{SBoCo}. However to simplify the
notation, we replace the number \(u_{pq}\) with the number
\(u_{qp}\) rather with the numbers \(\overline{u_{pq}}\). This
corresponds to that in \eqref{VGeE} we use the matrix \(U^{\ast}\)
rather than \(U\) as a matrix which parameterizes the set of all
\(\Omega_L\)-self-orthogonal subspaces.  The matrix \(U^\ast\) is
an arbitrary unitary matrix as well the matrix \(U\).
\end{remark}
\begin{definition}
\label{SAEMOp}%
Let \(U\) be a \(2\times2\) matrix. The operator \(\mathcal{L}_U\)
is
de\-fin\-ed  in the following way:\\
\hspace*{1.5ex}\textup{1}. %
\begin{minipage}[t]{0.92\linewidth}
The domain of definition \(\mathcal{D}_{\mathcal{L}_U}\) of the
operator \(\mathcal{L}_U\) is the set of all
\(x(t)\in\mathcal{D}_{\mathcal{L}_{\textup{max}}}\) which satisfy
the conditions \eqref{SBoCo1}-\eqref{SBoCo2}, \eqref{GBCo}.
\end{minipage}\\[1.2ex]
\hspace*{1.5ex}\textup{2}. %
\begin{minipage}[t]{0.92\linewidth}
For \(x\in\mathcal{D}_{\mathcal{L}_U}\), the action of the
operator \(\mathcal{L}_U\) is:
 \(\mathcal{L}_{U}x=\mathcal{L}_{\textup{max}}x\).
\end{minipage}
\end{definition}
\begin{remark}
\label{JDD}%
 In view of \eqref{SBoCo} and \eqref{SBoCo}, for any
matrix \(U\),
\[\mathcal{D}_{\mathcal{L}_{\textup{min}}}\subseteq\mathcal{D}_{\mathcal{L}_U}\,.\]
Thus for any matrix \(U\), the operator \(\mathcal{L}_U\) is an
extension of the operator \(\mathcal{L}_{\textup{min}}\):
\begin{equation}
\label{Betw}%
\mathcal{L}_{\textup{min}}\subseteq\mathcal{L}_U\subseteq
\mathcal{L}_{\textup{max}}\,.
\end{equation}
The equalities \eqref{SBoCo1} which determine the domain of
definition of the extension \(\mathcal{L}_U\) can be considered as
\emph{boundary conditions} posed on functions
\(x\in\mathcal{D}_{\mathcal{L}_{\textup{max}}}\). \textup{(}See
\textup{Remark \ref{GeBV}.)}
\end{remark}
The following Theorem is a reformulation of Lemma \ref{ODSS} in
the language of extensions of operators.
\begin{theorem}{\ } \\ %
\label{DSAE}
\hspace*{1.5ex}\textup{1}. %
\begin{minipage}[t]{0.92\linewidth}
If \(U\) is an unitary matrix, then the operator \(\mathcal{L}_U\)
is a selfadjoint differential operator which is an extension of the
minimal differential operator \(\mathcal{L}_{\textup{min}}\): %
\(\mathcal{L}_{\textup{min}}\subset\mathcal{L}_U\subset\mathcal{L}_{\textup{max}}\)\,.
\end{minipage}\\[1.0ex]
\hspace*{1.5ex}\textup{2}. %
\begin{minipage}[t]{0.92\linewidth}
Every differential operator \(\mathcal{L}\) which is a selfadjoint
extension of the minimal differential operator
\(\mathcal{L}_{\textup{min}}\), is of the form
\(\mathcal{L}=\mathcal{L}_U\) for some unitary matrix \(U\).\\
\end{minipage}
\hspace*{1.5ex}\textup{3}. %
\begin{minipage}[t]{0.92\linewidth}
For unitary matrices \(U_1,\,U_2\),
\begin{equation*}
(U_1=U_2)\Leftrightarrow(\mathcal{L}_{U_1}=\mathcal{L}_{U_2})\,.
\end{equation*}
\end{minipage}\\[2.0ex]
\end{theorem}
\noindent%
\textsf{Commutator of the operator \(\mathscr{F}_{[-a,a]}\) and \(\mathcal{L}_U\).}\\
 Let us
calculate the difference
\(\mathscr{F}_{[-a,a]}\mathcal{L}_{\textup{max}}x-%
\mathcal{L}_{\textup{max}}\mathscr{F}_{[-a,a]}x\) for
\(x\in\mathcal{D}_{\mathcal{L}_{\textup{max}}}\).
 Notice that
\(\mathcal{L}_{\textup{max}}x\in{}L^2([-a,a])\), so
\(\mathscr{F}_{[-a,a]}(\mathcal{L}_{\textup{max}}x)\) is defined.
Since \(x\in{}L^2([-a,a]\), the function
\(\mathscr{F}_{[-a,a]}x(t)\) is smooth on the closed interval
\([-a,a]\). (In fact this function is analytic in the whole real
axis.) All the more,
\(\mathscr{F}_{[-a,a]}x\in\mathcal{D}_{\mathcal{L}_{\textup{max}}}\).
Thus for \(x\in\mathcal{D}_{\mathcal{L}_{\textup{max}}}\), the
difference
\(\mathscr{F}_{[-a,a]}\mathcal{L}_{\textup{max}}x-%
\mathcal{L}_{\textup{max}}\mathscr{F}_{[-a,a]}x\)
is well defined.

Assuming that  \(x\in\mathcal{D}_{\mathcal{L}_{\textup{max}}}\) and that
\(-a<\alpha<\beta<a\),  we integrate
by parts twice\,%
\footnote{ {Like it is done in (2.31) of the manuscript
\cite{KaMa}.}
} :
\begin{multline}
\label{ffp}%
 \int\limits_{\alpha}^{\beta}
\left(-\frac{d\,\,}{d\xi}\Bigg(\bigg(1-\frac{\xi^2}{a^2}\bigg)\right)
\frac{dx(\xi)}{d\xi}\Bigg)e^{it\xi}d\xi=\\[1.0ex]
=-\bigg(1-\frac{\xi^2}{a^2}\bigg)\frac{dx(\xi)}{d\xi}\,e^{it\xi}%
\bigg|_{\xi=\alpha}^{\xi=\beta}
+it\bigg(1-\frac{\xi^2}{a^2}\bigg)x(\xi)%
e^{it\xi}\bigg|_{\xi=\alpha}^{\xi=\beta}-\\[1.0ex]
-it
\int\limits_{\alpha}^{\beta}x(\xi)\,\frac{d\,\,}{d\xi}
\left(\bigg(1-\frac{\xi^2}{a^2}\bigg)e^{it\xi}\right)\,d\xi\,.
\end{multline}

For \(x\in\mathcal{D}_{\mathcal{L}_{\textup{max}}}\), both limits
\(\displaystyle\lim_{t\to\pm{}a}(1-t^2/a^2)\frac{dx(t)}{dt}\) exist, are finite,
and
\begin{subequations}
\label{EVT}
\begin{align}
\label{EVT1}
\lim_{t\to{}-a}(1-t^2/a^2)\frac{dx(t)}{dt}&=\ \frac{2}{a^2}\,b_{-a}(x)\,,\\
\label{EVT2}
\lim_{t\to+a}(1-t^2/a^2)\frac{dx(t)}{dt}&=-\frac{2}{a^2}\,\,b_{\,a}(x)\,.\,
\end{align}
\end{subequations}
where \(b_{-a}(x), b_{a}(x)\) are defined in \eqref{GBCo} and also appear
in the boundary conditions \eqref{SBoCo}. Since the limits in \eqref{EVT}
are finite, we conclude that \(|x(t)|=O(\ln(a^2-t^2))\) as \(t\to\pm{}a,\,|t|<a\).
All the more, for \(x\in\mathcal{D}_{\mathcal{L}_{\textup{max}}}\)
\begin{equation}
\label{ZDC}
\lim_{t\to-a+0}\bigg(1-\frac{t^2}{a^2}\bigg)x(t)=0\,.
\end{equation}
Passing to the limit in \eqref{ffp} and taking into account \eqref{ZDC} and \eqref{EVT},
we obtain
\begin{multline}
\label{ffpOm}%
 \int\limits_{-a}^{a}
\left(-\frac{d\,\,}{d\xi}\Bigg(\bigg(1-\frac{\xi^2}{a^2}\bigg)
\frac{dx(\xi)}{d\xi}\Bigg)\right)e^{it\xi}d\xi=
\frac{2}{a}\bigg(b_{+}(x)e^{iat}+b_{-}(x)e^{-iat}\bigg)-
\\[1.0ex]
-it
\int\limits_{-a}^{a}x(\xi)\,\frac{d\,\,}{d\xi}
\left(\bigg(1-\frac{\xi^2}{a^2}\bigg)e^{it\xi}\right)\,d\xi\,.
\end{multline}
Transforming the last summand of the right hand side of \eqref{ffpOm},
we obtain
\begin{multline*}%
-it\int\limits_{-a}^{a}x(\xi)\,\frac{d\,\,}{d\xi}
\left(\bigg(1-\frac{\xi^2}{a^2}\bigg)e^{it\xi}\right)\,d\xi=\\[1.0ex]
=t^2\int\limits_{-a}^{a}x(\xi)e^{it\xi}\,d\xi+
\frac{it}{a^2}\int\limits_{-a}^{a}x(\xi)%
\frac{d\,\,}{d\xi}(\xi^2e^{it\xi})\,d\xi=
\end{multline*}%
\vspace{-2.0ex}
\hspace*{8.0ex}\bigg(\,\,since \(\dfrac{d\,\,}{d\xi}(\xi^2e^{it\xi})=\dfrac{d\,\,}{d\xi}%
\Big(-\dfrac{d^2\,\,}{dt^2}e^{it\xi}\Big)
=-\dfrac{d^2\,\,}{dt^2}\big(ite^{it\xi}\big)\)\,\,\bigg)\hspace*{2.0ex}
\vspace{2.0ex}
\begin{multline*}%
=t^2\int\limits_{-a}^{a}x(\xi)e^{it\xi}\,d\xi+\frac{t}{a^2}\,%
\frac{d^2\,\,}{dt^2}\bigg(t\int\limits_{-a}^{a}x(\xi)e^{it\xi}\,d\xi\bigg)=\\[1.0ex]
=t^2\int\limits_{-a}^{a}x(\xi)e^{it\xi}\,d\xi+\frac{d\,\,}{dt}%
\bigg(\frac{t^2}{a^2}\frac{d\,\,}{dt}\int\limits_{-a}^{a}x(\xi)e^{it\xi}\,d\xi\bigg)=
\end{multline*}
\begin{equation}
\label{Cont}
=t^2\int\limits_{-a}^{a}x(\xi)e^{it\xi}\,d\xi-
\frac{d\,\,}{dt}\Bigg(\bigg(1-\frac{t^2}{a^2}\bigg)\frac{d\,\,}{dt}%
\int\limits_{-a}^{a}x(\xi)e^{it\xi}\,d\xi\Bigg)-
\int\limits_{-a}^{a}\xi^2x(\xi)\,e^{it\xi}\,d\xi\,.
\end{equation}%
Unifying \eqref{ffpOm} and \eqref{Cont}, we obtain the equality
\begin{multline}%
\int\limits_{-a}^{a}
\left(\bigg(-\frac{d\,\,}{d\xi}\bigg(1-\frac{\xi^2}{a^2}\bigg)
\frac{d\,\,}{d\xi}+\xi^2\bigg)\,x(\xi)\right)e^{it\xi}\,d\xi=\\[1.0ex]
=\frac{2}{a}\bigg(b_{+}(x)e^{iat}+b_{-}(x)e^{-iat}\bigg)
+\Bigg(-\frac{d\,\,}{dt}\bigg(1-\frac{t^2}{a^2}\bigg)\frac{d\,\,}{dt}+t^2\Bigg)%
\int\limits_{-a}^{a}x(\xi)e^{it\xi}\,d\xi\,.
\end{multline}%
We summarize the above calculation as
\begin{lemma}
Let \(\mathscr{F}_{[-a,a]}\) be the Fourier operator truncated on
the finite symmetric interval \([-a,a]\). Let
\(\mathcal{L}_{\textup{max}}\) be
the maximal differential operator with domain of definition %
\(\mathcal{D}_{\mathcal{L}_{\textup{max}}}\)
generated by the formal differential operator
\(\displaystyle{}L=-\frac{d\,\,}{dt}\bigg(1-\frac{t^2}{a^2}\bigg)\frac{d\,\,}{dt}+t^2\).
\textup{(}See \textup{Definition \ref{MaxDO}.)}

If  \(x\in\mathcal{D}_{\mathcal{L}_{\textup{max}}}\),
 then  \(\mathscr{F}_{[-a,a]}x\in\mathcal{D}_{\mathcal{L}_{\textup{max}}}\),
and the equality holds
\begin{equation}
\label{ComRel}
(\mathscr{F}_{[-a,a]}\mathcal{L}_{\textup{max}}x)(t)-
(\mathcal{L}_{\textup{max}}\mathscr{F}_{[-a,a]}x)(t)=
\frac{2}{a}\bigg(b_{+}(x)e^{iat}+b_{-}(x)e^{-iat}\bigg)\,.
\end{equation}
\end{lemma}
Every selfadjoint differential operator generated by the formal
differential operator \(L\) is a restriction of the maximal
differential operator \(\mathcal{L}_{\textup{max}}\) on the
appropriate domain of definition. According to Theorem~\ref{DSAE},
the set of such self-adjoint operators coincides with the set of
operators \(\mathcal{L}_U\), where \(U\) is an arbitrary
\(2\times2\) unitary matrix. The domain of definition
\(\mathcal{D}_{\mathcal{L}_U}\) of the operator \(\mathcal{L}_U\)
is distinguished from the domain
\(\mathcal{D}_{\mathcal{L}_{\textup{max}}}\) by the boundary
conditions \eqref{SBoCo} constructed from \(U\). The next theorem
answers the question which operators \(\mathcal{L}_U\) commute
with the truncated Fourier operator \(\mathscr{F}_{[-a,a]}\).
\begin{theorem}{\ }\\[0.5ex]
\label{WECF}
\hspace*{1.0ex}\textup{1}. %
\begin{minipage}[t]{0.94\linewidth}
If \(U=I\), where \(I\) is \(2\times2\)  identity matrix, then
the differential operator\,%
\footnotemark\,%
\(\mathcal{L}_I\) commutes with the truncated\,\footnotemark\,%
 Fourier operator~%
\(\mathscr{F}_{[-a,a]}\):
\begin{equation}
\label{DCoRe}
\mathscr{F}_{[-a,a]}\mathcal{L}_I\,x=\mathcal{L}_I\mathscr{F}_{[-a,a]}\,x\quad \forall\, x\in%
\mathcal{D}_{\mathcal{L}_I}\,.
\end{equation}
\end{minipage}\\[1.5ex]
\hspace*{1.0ex}\textup{2}. %
\begin{minipage}[t]{0.94\linewidth}
If \,\,\(U\not=I\), then the operator \(\mathcal{L}_U\) do not commute with the operator
\(\mathscr{F}_{[-a,a]}\):\\[1.0ex]
\hspace{1.0ex}%
\hspace{1.0ex}%
\textup{a.}
\begin{minipage}[t]{0.93\linewidth}
There exist vectors \(x\in\mathcal{D}_{\mathcal{L}_U}\) such that
\(\mathscr{F}_{[-a,a]}\in\mathcal{D}_{\mathcal{L}_U}\), so both operators
\(\mathscr{F}_{[-a,a]}\mathcal{L}_{U}\) and \(\mathcal{L}_{U}\mathscr{F}_{[-a,a]}\)
are applicable to \(x\), but
\begin{equation}
\mathscr{F}_{[-a,a]}\mathcal{L}_{U}x\not=
\mathcal{L}_{U}\mathscr{F}_{[-a,a]}x\,;
\end{equation}
\end{minipage}\\[1.0ex]
\textup{b.}
\begin{minipage}[t]{0.93\linewidth}
 There exist vectors \(x\in\mathcal{D}_{\mathcal{L}_U}\) such that
\(\mathscr{F}_{[-a,a]}x\not\in\mathcal{D}_{\mathcal{L}_U}\), so the operator
\(\mathcal{L}_{U}\mathscr{F}_{[-a,a]}\) even can not be applied to such \(x\).
\end{minipage}\\
\end{minipage}\\
\addtocounter{footnote}{-1}
\footnotetext{\,\(\mathcal{L}_I=\mathcal{L}_U \textup{ for }U=I\).}
\addtocounter{footnote}{1}
\footnotetext{\,\(\mathscr{F}_{[-a,a]}=F_{E}\) for \(E=[-a,a]\).}%
\end{theorem}
\begin{proof}{\ }\\
\textbf{1}.
For \(U=I\), the boundary conditions \eqref{SBoCo} take the form
\begin{equation}
\label{DBoCo}
b_{-a}(x)=0,\quad b_{a}(x)=0\,.
\end{equation}
Thus, the domain of definition \(\mathcal{D}_{\mathcal{L}_I}\) of the operator
\(\mathcal{L}_I\) is:
\begin{equation}
\label{DDLI}
\mathcal{D}_{\mathcal{L}_I}=\big\lbrace{}x:\,x\in\mathcal{D}_{\mathcal{L}_{\textup{max}}},\
b_{-a}(x)=0,\,b_{a}(x)=0\big\rbrace\,.
\end{equation}
Every smooth function \(x(t)\) on \((-a,a)\) which derivative is
bounded: \(\sup_{t\in(-a,a)}|x^\prime(t)|<\infty\), belongs to
\(\mathcal{D}_{\mathcal{L}_{\textup{max}}}\). Moreover, according
to \eqref{GBCo}, every such a function  satisfies the boundary
condition \eqref{DDLI}, i.e. \(b_{-a}(x)=0\),\,\(b_{a}(x)=0\).
Hence \emph{every smooth function on \((-a,a)\) which derivative
is bounded on \((-a,a)\), belongs to domain of definition
\(\mathcal{D}_{\mathcal{L}_I}\) of the operator
\(\mathcal{L}_I\)}. In particular, if \(x\in{}L^2((-a,a)\) and
\(y=\mathscr{F}_{[-a,a]}x\), then
\(y\in\mathcal{D}_{\mathcal{L}_I}\).
Thus for \(x\in\mathcal{D}_{\mathcal{L}_I}\) both summands in the expression %
\(\mathscr{F}_{[-a,a]}\mathcal{L}_Ix-\mathcal{L}_I\mathscr{F}_{[-a,a]}x\)
are well defined. Since the operator \(\mathcal{L}_I\) is a restriction of the
operator \(\mathcal{L}_{\textup{max}}\), then
\begin{equation*}%
\mathscr{F}_{[-a,a]}\mathcal{L}_Ix-\mathcal{L}_I\mathscr{F}_{[-a,a]}x=
\mathscr{F}_{[-a,a]}\mathcal{L}_{\textup{max}}x-%
\mathcal{L}_{\textup{max}}\mathscr{F}_{[-a,a]}x\,\,\textup{ for }
x\in\mathcal{D}_{\mathcal{L}_I}\,.
\end{equation*}
In view of \eqref{ComRel} and \eqref{DBoCo}, the equality \eqref{DCoRe} holds.\\[1.0ex]
\textbf{2}. Let \(U\not=I\). Then at least of one value \(u_{11}-1\) or \(u_{22}-1\)
differs from zero. For definiteness, let \ \(u_{11}-1\not=0\). Set
\begin{equation}%
\label{SpX}%
\gamma=\frac{1+u_{11}}{i(1-u_{11})}\,,\quad x(t)=\psi_{-}(t)+\gamma\varphi_{-}(t)+x_0(t),
\end{equation}
where \(x_0(t)\) is a smooth function which support is a compact subset of the
\emph{open} interval \((-a,a)\):
\begin{equation}%
\label{SuppX}%
\textup{supp}\,x_0\Subset(-a,a)\,.
\end{equation}
The function \(x_0\) will be chosen later. According to \eqref{Rewr}, \eqref{SuppX}
and the choice of \(\gamma\), \emph{for any choice of} \(x_0(t)\),
the function \(x(t)\) from \eqref{SpX} satisfy the
boundary conditions \eqref{SBoCo}.
Thus,
\begin{equation}
x(t)\in\mathcal{D}_{\mathcal{L}_U}\,.
\end{equation}
for any choice of \(x_0\). Moreover
\begin{equation}
\label{but}
b_{-a}(x)=1,\quad  b_{a}(x)=0\,.
\end{equation}
For the function \(y(t)=(\mathscr{F}_{(-a,a)}x)(t)\), \emph{the
boundary conditions \eqref{SBoCo} either hold, or does not hold.}
This depends on the choice of the function \(x_0\). If
\eqref{SBoCo} hold for this \(y\), then
\(\mathscr{F}_{(-a,a)}x\in\mathcal{D}_{\mathcal{L}_{U}}\) and the
equality \eqref{ComRel} can be interpreted as the equality
\begin{equation}
\label{Exam1}
(\mathscr{F}_{(-a,a)}\mathcal{L}_{U}x)(t)-
(\mathcal{L}_{U}\mathscr{F}_{(-a,a)}x)(t)=
\frac{2}{a}\bigg(b_{+}(x)e^{iat}+b_{-}(x)e^{-iat}\bigg)\,.
\end{equation}
In view of \eqref{but}, \((\mathscr{F}_{(-a,a)}\mathcal{L}_{U}x)(t)-
(\mathcal{L}_{U}\mathscr{F}_{(-a,a)}x)(t)\not=0\).

Let us show that both of the possibilities %
\(\mathscr{F}_{(-a,a)}x\in\mathcal{D}_{\mathcal{L}_{U}}\) and
\(\mathscr{F}_{(-a,a)}x\not\in\mathcal{D}_{\mathcal{L}_{U}}\)
are realizable. Since the function \(y(t)=(\mathscr{F}_{(-a,a)}x)(t)\)
is smooth on \([-a,a]\),
\begin{equation*}%
b_{-a}(y)=0,\,b_{a}(y)=0,\,\,c_{-a}(y)=-y(-a),\,c_{-a}(y)=-y(a)\,.
\end{equation*}%
Thus as applied to the function \(y\), the boundary conditions \eqref{SBoCo}
take the form
\begin{subequations}
\label{SpBC}
\begin{gather}
\label{SpBC1}
(1-u_{11})y(-a)-u_{12}y(a)=0\,,\\
\label{SpBC2}
u_{21}y(-a)-(1-u_{22})y(a)=0\,.
\end{gather}
\end{subequations}
If, using the freedom of choice of the function \(x_0(t)\) in \eqref{SpX}, we can arbitrary
prescribe the values \(y(-a)\) and \(y(a)\), then we can either satisfy the
boundary conditions \eqref{SpBC} (prescribing \(y(-a)=0,\,y(a)=0\)),
or violate them (if \(u_{11}\not=1\),
we prescribe \(y(-a)=1,\,y(a)=0\), if  \(u_{11}\not=1\),
we prescribe \(y(-a)=0,\,y(a)=1\).) The reference to Lemma below
finishes the proof.
\end{proof}
\begin{lemma}
Given complex numbers \(y_1\) and \(y_2\), there exists a smooth
function
\(x_0(t)\) on \([-a,a]\) which possesses the properties:\\[1.0ex]
\hspace*{1.5ex}\textup{1}. \(\hspace*{20.0ex}\textup{supp}\,x_0\Subset(-a,a)\,.\)\\[1.0ex]
\hspace*{1.5ex}\textup{2}. \(\hspace*{10.0ex}y_0(-a)=y_1,\ y_0(a)=y_2\), where \(y_0=\mathscr{F}_{[-a,a]}(x_0)\).
\end{lemma}
\begin{proof}
The evaluations \(y(-a)\) and \(y(-a)\) are linearly independent linear functionals
on the space of functions on \((-a,a)\) which are smooth and compactly supported:
\[y(-a)=\frac{1}{\sqrt{2\pi}}\int\limits_{-a}^{a}x(\xi)e^{-ia\xi}\,d\xi,\quad
y(a)=\frac{1}{\sqrt{2\pi}}\int\limits_{-a}^{a}x(\xi)e^{ia\xi}\,d\xi\,,
\]
and the functions \(e^{-ia\xi},\,e^{ia\xi}\) generating these linear functionals
are linearly independent on any non-empty open subinterval of the interval \((-a,a)\).
\end{proof}{\ }\\
\noindent
\textsf{Properties of the operator \(\mathcal{L}_I\).} As we have established,
Theorems \ref{DSAE} and \ref{WECF}, the only selfadjoint differential
operator which is generated by the formal operator \(L\) and which commutes
with the truncated Fourier operator \(\mathscr{F}_{[-a,a]}\) is the operator %
\(\mathcal{L}_I\). Let as discuss properties of the operator \(\mathcal{L}_I\).

The following lemma gives an alternative definition of the domain~%
\(\mathcal{D}_{\mathcal{L}_{I}}\).
\begin{lemma}
\label{CAEP}%
 Let a function \(x(t)\) belong to
 \(\mathcal{D}_{\mathcal{L}_{\textup{max}}}\).
 Then the function \(x(t)\) belong to
 \(\mathcal{D}_{\mathcal{L}_{I}}\) if and only if both
limits
\begin{equation}
\label{EPV}
x(-a)=\lim_{t\to{-a+0}}x(t),\quad x(a)=\lim_{t\to{a-0}}x(t)\,.
\end{equation}
exist and are finite.
\end{lemma}
\begin{proof} \ \textsf{1}.\,\,%
Functions \(x(t)\) belonging to \(\mathcal{D}_{\mathcal{L}_I}\) possess the properties
\begin{align}
\label{Pr1}
\lim_{t\to\pm{a}}(a^2-t^2)\frac{dx(t)}{dt}&=0,\\
\label{Pr2}
\int\limits_{-a}^{a}\bigg|\frac{d\,\,}%
{d\xi}\bigg((a^2-\xi^2)\frac{dx(\xi)}{d\xi}\bigg)\bigg|^2\,&d\xi=C^2<\infty, \ \ C>0.
\end{align}
From \eqref{Pr2} and the Schwarz inequality it follows that
\begin{equation*}%
\int\limits_{t_1}^{t_2}\bigg|\frac{d\,\,}%
{d\xi}\bigg((a^2-\xi^2)\frac{dx(\xi)}{d\xi}\bigg)\bigg|\,d\xi\leq{}C\,\sqrt{t_2-t_1}\,,
\quad -a<t_1<t_2<a\,.
\end{equation*}%
All the more,
\[\Bigg|(a^2-\xi^2)\frac{dx(\xi)}{d\xi}\bigg|_{\xi=t_1}^{\xi=t_2}\Bigg|\leq
C\,\sqrt{t_2-t_1}\,.\]
We use the last inequality for \(t_1=-a+0,\,t_2=t\), where \mbox{\(-a<t<a\).}
Taking into account \eqref{Pr1}, we deduce that
\[\bigg|(a^2-t^2)\frac{dx(t)}{dt}\bigg|\leq{}C\sqrt{a+t}\,,\quad -a\leq{}t\leq{}0\,.\]
Analogously,
\[\bigg|(a^2-t^2)\frac{dx(t)}{dt}\bigg|\leq{}C\sqrt{a-t}\,,\quad 0\leq{}t\leq{}a\,.\]
From two last inequalities it follows that
\[\bigg|(a^2-t^2)\frac{dx(t)}{dt}\bigg|\leq\frac{1}{\sqrt{a}}\sqrt{a^2-t^2},\quad\,-a<t<a\,.\]
Finally, from \eqref{Pr1} and \eqref{Pr2} we deduced the inequalities
\begin{equation}
\bigg|\frac{dx(t)}{dt}\bigg|\leq{}\frac{C}{\sqrt{a}}\frac{1}{\sqrt{a^2-t^2}}\,,\quad -a<t<a\,.
\end{equation}
and
\begin{equation}
\label{MoCo}
|x(t_2)-x(t_1)|\leq{}\frac{C}{\sqrt{a}} \int\limits_{t_1}^{t_2}\frac{d\xi}{\sqrt{a^2-\xi^2}},\quad
-a<t_1<t_2<a\,.
\end{equation}
Since \(\int\limits_{-a}^{a}\frac{dt}{\sqrt{a^2-t^2}}<\infty\),
the function \(x(t)\) is uniformly continuous on the interval
\((-a,a)\). Therefore the limits \eqref{EPV} exist end are
finite.\\
\hspace*{1.5ex}\textsf{2}. According to Lemma \ref{GBCL}, both
limits \(\lim_{t\to{\pm{}(a-0)}}(t\pm{}a)\frac{dx(t)}{dt}\) exist
and are finite. If \(x\not\in\mathcal{D}_{\mathcal{L}_I}\), them
at least one of these limits is not zero. If for example
\(\lim_{t\to{}a-0}(t-a)\frac{dx(t)}{dt}\not=0\), then the function
\(x(t)\) grows logarithmically as \(t\to{a-0}\).
\end{proof}
\begin{theorem}{\ }\\[0.5ex]
\label{DiSp}%
\hspace*{1.5ex}\textup{1}.
\begin{minipage}[t]{0.95\linewidth}
The selfadjoint operator \(\mathcal{L}_I\) is an operator with  discrete spectrum.
\end{minipage}\\[1.0ex]
\hspace*{1.5ex}\textup{2}.
\begin{minipage}[t]{0.95\linewidth}
The spectrum of the operator \(\mathcal{L}_I\) is a sequence of positive eigenvalues
of multiplicity one which tends to \(+\infty\).
\end{minipage}\\[1.0ex]
\hspace*{1.5ex}\textup{3}.
\begin{minipage}[t]{0.95\linewidth}
The number \(\lambda\) is an eigenvalue of the differential
operator \(\mathcal{L}_I\) if and only if there exists the
non-zero solution \(e(t,\lambda)\) of the boundary value problem
for the differential equation
\begin{subequations}
\label{EiVaPr}
\begin{equation}%
\label{EiVaPr1}%
-\frac{d}{dt}\Bigg(\bigg(1-\frac{t^2}{a^2}\bigg)\frac{de(t,\lambda)}{dt}\Bigg)+
t^2e(t,\lambda)=\lambda{}e(t,\lambda)
\end{equation}
with the boundary conditions\\[-1.5ex]
\begin{equation}%
\label{EiVaPr2}%
 e(-a,\lambda)\textup{ is finite }, \ \
e(a,\lambda)\textup{ is finite. }
\end{equation}%
\end{subequations}
This solution \(e(t,\lambda)\) is an eigenvector of the operator
\(\mathcal{L}_I\) corresponding to the eigenvalue \(\lambda\).
\end{minipage}
\end{theorem}
\begin{remark}
\label{FPrSph}%
The solutions of the boundary value problem \eqref{EiVaPr} are
known as the \emph{prolate spheroidal wave functions}. There is a
literature where these functions are discussed and studied. See
for example \textup{\cite{ChSt}, \cite{Fl}, \cite{KPS},
\cite{MSch}, \cite{SMCLC}.}
\end{remark}
 If \(A\) is a symmetric operator in a Hilbert space \(\mathfrak{H}\)
which domain of definition \(\mathcal{D}_A\) is dense in \(\mathfrak{H}\)
and \(M\) is a bounded selfadjoint operator defined everywhere in \(\mathfrak{H}\)\,,
then the  operators \(A\) and \(B=A+M\) (\(\mathcal{D}_{B}=\mathcal{D}_{A}\)) are selfadjoint
or not simultaneously, and  spectra of \(A\) and \(B\) are discrete or not
simultaneously.

We use this fact in the case when \(\mathfrak{H}=L^2((-a,a))\), \(A=\mathcal{L}_I\),
\(Mx(t)=x(t)-t^2x(t)\), so the operator \(B\) is a differential operator \(\Lambda\)
of the form
\begin{equation}
\label{Leg}
(\Lambda{}x)(t)=-\frac{d\,\,}{dt}\bigg(\Big(1-\frac{t^2}{a^2}\Big)\frac{dx(t)}{dt}\bigg)+x(t)\,.
\end{equation}
which domain of definition \(\mathcal{D}_{\Lambda}\) coincides with the domain of definition
\(\mathcal{D}_{\mathcal{D}_{\mathcal{L}_I}}\) of the operator \(\mathcal{L}_I\).
(See \eqref{DDLI} and \eqref{maxdo1}.)
\begin{lemma}
\label{Nonne}%
Each of the operators \(\mathcal{L}_I\) and \(\Lambda\) is non-negative, and
for every \(x\in\mathcal{D}_{\mathcal{L}_I}=\mathcal{D}_{\Lambda}\) the equalities hold:
\begin{gather}
\label{none1}
\langle\mathcal{L}_Ix,x\rangle=\int\limits_{-a}^{a}\bigg(1-\frac{\xi^2}{a^2}\bigg)%
\bigg|\frac{dx(\xi)}{d\xi}\bigg|^2d\xi+\int\limits_{-a}^{a}\xi^2|x(\xi)|^2d\xi\,,\\
\label{none2}
\langle\Lambda{}x,x\rangle=\int\limits_{-a}^{a}\bigg(1-\frac{\xi^2}{a^2}\bigg)%
\bigg|\frac{dx(\xi)}{d\xi}\bigg|^2d\xi+\int\limits_{-a}^{a}|x(\xi)|^2d\xi\,.
\end{gather}
\end{lemma}
\begin{proof} Let \(-a<\alpha<\beta<a\).
Integrating by parts we obtain
\begin{multline*}%
\int\limits_{\alpha}^{\beta}\bigg(-\frac{d\,\,}{d\xi}%
\bigg(1-\frac{\xi^2}{a^2}\bigg)\frac{dx(\xi)}{d\xi}\bigg)\overline{x(\xi)}\,d\xi=\\[1.0ex]
=-\bigg(1-\frac{\xi^2}{a^2}\bigg)\frac{dx(\xi)}{d\xi}\cdot\overline{x(\xi)}\bigg|_{\xi=\alpha}^{\xi=\beta}
+\int\limits_{\alpha}^{\beta}\bigg(1-\frac{\xi^2}{a^2}\bigg)\frac{dx(\xi)}{d\xi}\cdot
\overline{\frac{dx(\xi)}{d\xi}}\,\,d\xi\,.
\end{multline*}
According to the boundary conditions \eqref{DBoCo}, \[\lim_{t\to\pm{}(a-0)}
\bigg(1-\frac{t^2}{a^2}\bigg)\frac{dx(t)}{dt}=0\,,\]
According to Lemma \ref{CAEP}, %
\[|x(t)|=O(1) \ \textup{ as } \ |t|\to{}a-0\,.\]
Passing to the limit as \(\alpha\to-a+0,\,\beta\to{}a-0\), we obtain the equality
\begin{multline}
\label{QFo}%
\int\limits_{-a}^{a}\bigg(-\frac{d\,\,}{d\xi}%
\bigg(1-\frac{\xi^2}{a^2}\bigg)\frac{dx(\xi)}{d\xi}\bigg)\overline{x(\xi)}\,d\xi=
\int\limits_{-a}^{a}\bigg(1-\frac{\xi^2}{a^2}\bigg)\bigg|\frac{dx(\xi)}{d\xi}\bigg|^2\,d\xi\,,\\
\textup{for every }\ x\in\mathcal{D}_{\mathcal{L}_I}=\mathcal{D}_{\Lambda}\,.
\end{multline}
\end{proof}
\begin{proof}[Proof of \textup{Theorem \ref{DiSp}}.] Let
\begin{equation}
\mathfrak{B}=\Big\lbrace{}x\in\mathcal{D}_{\Lambda}:\,%
\langle\Lambda{}x,\Lambda{}x\rangle_{L^2(-a,a)}\leq{}1\Big\rbrace\,.
\end{equation}
be a preimage of the unit ball of the space \(L^2(-a,a)\) with
respect to the mapping \(x\to\Lambda{}x\). To prove that the
spectrum of \(\Lambda\) is discrete it is enough to prove that the
set \(\mathfrak{B}\) is precompact in \(L^2(-a,a)\).
The condition \(\langle\Lambda{}x,\Lambda{}x\rangle_{L^2(-a,a)}\leq{}1\)
for a function \(x\in\mathcal{D}_{\Lambda}\) means that
\begin{equation}%
\label{prei}
\int_{-a}^{a}\bigg|-\frac{d\,}{d\xi}\bigg(\Big(1-\frac{\xi^2}{a^2}\Big)
\frac{dx(\xi)}{d\xi}\bigg)+x(\xi)\bigg|^2\,d\xi\leq1
\end{equation}
In view of \eqref{QFo},
\begin{multline*}%
\int_{-a}^{a}\bigg|-\frac{d\,}{d\xi}\bigg(\Big(1-\frac{\xi^2}{a^2}\Big)
\frac{dx(\xi)}{d\xi}\bigg)\bigg|^2\,d\xi+\int\limits_{-a}^{a}|x(\xi)|^2\,d\xi\leq\\
\int_{-a}^{a}\bigg|-\frac{d\,}{d\xi}\bigg(\Big(1-\frac{\xi^2}{a^2}\Big)
\frac{dx(\xi)}{d\xi}\bigg)+x(\xi)\bigg|^2\,d\xi
\,.
\end{multline*}%
Therefore from \eqref{prei} it follows that
\begin{gather}
\label{UnBo}
\int\limits_{-a}^{a}|x(\xi)|^2\,d\xi\leq1\\
\intertext{and}
\label{EqCon}
\int_{-a}^{a}\bigg|-\frac{d\,}{d\xi}\bigg(\Big(1-\frac{\xi^2}{a^2}\Big)
\frac{dx(\xi)}{d\xi}\bigg)\bigg|^2\,d\xi\leq{}1\,.
\end{gather}
Inequality \eqref{EqCon} the inequality \eqref{Pr2} for \(C=a^2\).
According to \eqref{MoCo}, the function \(x\) satisfy the inequality
\begin{equation}%
\label{UBE}
|x(t_2)-x(t_1)|\leq{}a^{3/2}
\int\limits_{t_1}^{t_2}\frac{d\xi}{\sqrt{a^2-\xi^2}},\quad
-a<t_1<t_2<a\,.
\end{equation}
Thus, the set of the functions \(x\) belonging to \(\mathfrak{B}\) is uniformly
bounded, \eqref{UnBo}, and equicontinuous, \eqref{UBE}. Therefore, the set
\(\mathfrak{B}\) is precompact in \(L^2([-a,a])\).

Thus the spectra of the operators \(\Lambda\) and \(\mathcal{L}_I\) is discrete, i.e. consists of
isolated eigenvalues. According to \eqref{none1}, the eigenvalues of the
operator \(\mathcal{L}_I\) are positive. If \(\lambda\) is an eigenvalue
of the operator \(\mathcal{L}_I\) and \(e(t,\lambda)\) is an eigenfunction
which corresponds to this \(\lambda\), then, since
\(e(t,\lambda)\in\mathcal{D}_{\mathcal{L}_I}\), the function \(e(t,\lambda)\)
is continuous in \(t\) at the points \(t=a\) and \(t=-a\).
(Lemma \ref{CAEP}.)

Moreover the function
\(e(t,\lambda)\) is the solution of the differential equation \(Lx=\lambda{}x\).
As any  solution of this equation, the function \(e(t,\lambda)\) is a linear combination
of the solutions  \(x_{1}^{-}(t,\lambda)\) and
\(x_{2}^{-}(t,\lambda)\). (The solutions \(x_{1}^{\pm}(t,\lambda), x_{2}^{\pm}(t,\lambda)\)
were introduced in Lemma \ref{ABSNS}).
From the behavior of the functions
\(e(t,\lambda),\,x_{1}^{-}(t,\lambda),\,x_{2}^{-}(t,\lambda)\) by \(t\to{}-a+0\)
we deduce that the function \(e(t,\lambda)\)  is proportional to \(x_{1}^{-}(t,\lambda)\):
\[e(t,\lambda)=C_{-}x_{1}^{-}(t,\lambda),\,\quad C_{-}\not=0\ \ \textup{is a constant}\,.\]
Analogously,
\[e(t,\lambda)=C_{+}x_{1}^{+}(t,\lambda),\,\quad C_{+}\not=0\ \ \textup{is a constant}\,.\]
Thus, up to the proportionality, there is only one eigenfunction corresponding
to the eigenvalue \(\lambda\).
\end{proof}
\begin{remark}
\label{ChE}
Thus if \(\lambda\) is an eigenvalue of \(\mathcal{L}_I\), then
\(C_{-}x_{1}^{-}(t,\lambda)=C_{+}x_{1}^{+}(t,\lambda)\).
Since the differential equation \(Lx=\lambda{}x\) is invariant with respect to
the change of variable \(t\to-t\), then the functions \(e(-t,\lambda)\),
\(e(t,\lambda)\pm{}e(-t,\lambda)\) are eigenfunctions as well.
Since there in only on eigenfunction up to proportionality,
then either \(e(t,\lambda)=e(-t,\lambda)\), or \(e(t,\lambda)=-e(-t,\lambda)\).
Thus, either \(C_{+}=C_{-}\), or \(C_{+}=-C_{-}\).
\end{remark}
\begin{remark}
\label{ChEr} The spectral analysis of the operator \(\Lambda\) can
be done explicitly. Its eigenfunctions are essentially the
Legandre polynomials, the spectrum also can by found  explicitely.
The property of the spectrum of \(\Lambda\) to be discrete may be
derived from this analysis. However we prefer to present less
explicit but more general reasoning.
\end{remark}

\vspace{4.0ex}
\begin{minipage}[h]{0.45\linewidth}
Victor Katsnelson\\[0.2ex]
Department of Mathematics\\
The Weizmann Institute\\
Rehovot, 76100, Israel\\[0.1ex]
e-mail:\\
{\small\texttt{victor.katsnelson@weizmann.ac.il}}
\end{minipage}

\vspace{3.0ex}
\begin{minipage}[h]{0.45\linewidth}
Ronny Machluf\\[0.2ex]
Department of Mathematics\\
The Weizmann Institute\\
Rehovot, 76100, Israel\\[0.1ex]
e-mail:\\
\texttt{ronny-haim.machluf@weizmann.ac.il}
\end{minipage}
\end{document}